%% file: main.tex
\pdfoutput=1
\documentclass[11pt,nodots,reqno,eqnsec]{style/orca}

\usepackage{bigints}
\usepackage{bookmark}
\usepackage{caption}
\usepackage{graphicx}
\usepackage{subcaption}
\usetikzlibrary{arrows.meta, decorations.markings, calc, positioning, decorations.pathmorphing, decorations.shapes, shapes.geometric, calligraphy, patterns}
\tikzset{->-/.style={decoration={
  markings,
  mark=at position .5 with {\arrow{>}}},postaction={decorate}}}
\tikzset{-<-/.style={decoration={
  markings,
  mark=at position .5 with {\arrow{<}}},postaction={decorate}}}

\DeclareSymbolFont{Letters}{OML}{cmm}{m}{it}
\DeclareMathSymbol{\psi}{\mathalpha}{Letters}{"20}
\newtheoremstyle{dotless}{}{}{\itshape}{}{\bfseries}{}{ }{}
\theoremstyle{dotless}

\newcounter{tpointn}[section]
\renewcommand{\thetpointn}{%
  \ifnum\value{section}=0 %
    \else\thesection.\arabic{tpointn}\fi
  }
\renewcommand{\tpointn}[2][]{\vspace{2mm}\par%
  \ifx&#1&\refstepcounter{tpointn}\fi%
  \noindent{\bf #2\ifx&#1& \thetpointn\fi\if@dots.\fi\, ---}}

\renewcommand{\thesubsubsection}{\thesubsection.\arabic{subsubsection}}

\graphicspath{{figures}}
\input{figures/figures}

\bibliography{main}

\title{Completing the \lowercase{\texorpdfstring{$c_2$}{c2}} completion conjecture for \lowercase{\texorpdfstring{$p=2$}{p=2}}}
\date{\today}
\author{Simone Hu}
\address{Mathematical Institute \newline
  \indent University of Oxford \newline
  \indent Oxford, United Kingdom, OX2 6GG}
\email{simone.hu@maths.ox.ac.uk}

\author{Karen Yeats}
\address{Department of Combinatorics and Optimization \newline
  \indent Faculty of Mathematics, University of Waterloo \newline
  \indent Waterloo, ON, Canada, N2L 3G1}
\email{kayeats@uwaterloo.ca}

\begin{document}

\begin{abstract}
  The $c_2$-invariant is an arithmetic graph invariant useful for understanding Feynman periods.
Brown and Schnetz conjectured that the $c_2$-invariant has a particular symmetry known as completion invariance. This paper will prove completion invariance of the $c_2$-invariant in the $p=2$ case, extending previous work of one of us.
The methods are combinatorial and enumerative involving counting certain partitions of the edges of the graph.
\end{abstract}

\maketitle

\setcounter{tocdepth}{1}
\tableofcontents
\setcounter{tocdepth}{3}

\input{intro}
\input{background}

\input{t-case}

\input{s-case}

\input{r-case}

\input{completing}
\input{higherp}

\printbibliography
\bigskip
\setlength\parskip{0pt}

\end{document}

%% file: figures/figures.tex
\newcommand{\Blob}{%
  \begin{tikzpicture}[main node/.style={inner sep=0,minimum size =.08cm,circle,fill=black!80,draw},node distance = 0.75cm,scale=0.8]
    \node[main node, label={[shift={(0,-2em)}]\scriptsize$u_2$}] (f) at (0,0) {};
    \node[main node, label=left:\scriptsize$u_1$] (g) at (-1,0) {};
    \node[main node, label=right:\scriptsize$u_3$] (i) at (1,0) {};
    \node[style={inner sep=2,circle,draw}] (v) at (0,1) {};
    \node[main node, label=above:\scriptsize$u$] (v) at (0,1) {};

    \draw (v) -- node[below left=0.01pt] {\scriptsize $2$} (f);
    \draw (v) -- node[left=0.1pt]{\scriptsize $1$} (g);
    \draw (v) -- node[right=0.1pt]{\scriptsize $3$} (i);
    \draw (g) to [bend right=60] (f) to [bend right=60] (i);
    \draw (g) to [bend right=110, in=290, looseness=2.0] (i);
  \end{tikzpicture}
}

\newcommand{\BlobBase}{%
  \begin{tikzpicture}[main node/.style={inner sep=0,minimum size =.08cm,circle,fill=black!80,draw},node distance = 0.75cm,scale=0.8]
    \node[main node, label=above:\scriptsize$u_2$] (f) at (0,0) {};
    \node[main node, label=above:\scriptsize$u_1$] (g) at (-1,0) {};
    \node[main node, label=above:\scriptsize$u_3$] (i) at (1,0) {};

    \draw (g) to [bend right=60] (f) to [bend right=60] (i);
    \draw (g) to [bend right=110, in=290, looseness=2.0] (i);
  \end{tikzpicture}
}

\newcommand{\BlobTwo}{%
  \node[main node, label={[shift={(0,-1.75em)}]\scriptsize $u_2$}] (f) at (0,-0.30) {};
  \node[main node, label=above:{\scriptsize$u=u_1 = u_3$}] (g) at (0,0.25) {};

  \draw (g) to [out=195, in=90, looseness=1] (-0.75,-0.25) to [out=270, in=225, looseness=1.5] (f);
  \draw (f) to [out=315, in=270, looseness=1.5] (0.75,-0.25) to [out=90, in=345, looseness=1] (g);
  \draw (g) to [out=188, in=180, looseness=2.25] (0,-1.25) to [out=0, in=352, looseness=2.25] (g);
}

\newcommand{\BlobThree}{%
  \node[main node, label=above:{\scriptsize$u_1$}] (i) at (-0.5,0) {};
  \node[main node, label=above:{\scriptsize$u=u_2=u_3$}] (h) at (1,0) {};

  \draw (h) to [out=225, in=320, looseness=25] (h);
  \draw (h) to [bend left=45] (i);
  \draw (h) to [out=330, in=0, looseness=2] (0.25,-1.25);
  \draw (0.25,-1.25) to [out=180, in=230, looseness=1.5] (i);
}

\newcommand{\ExampleBlob}{%
  \begin{tikzpicture}[main node/.style={inner sep=0,minimum size =.08cm,circle,fill=black!80,draw},node distance = 0.75cm,scale=0.8]
    \clip (-1.5,-1.5) rectangle (8,0.75);
    \begin{scope}
      \BlobThree
    \end{scope}
    \node at (2.75,-0.5) {$\bigcap$};
    \begin{scope}[shift={(5,0)}]
      \BlobTwo
    \end{scope}
  \end{tikzpicture}
}

\newcommand{\ShapeOne}{%
  \begin{tikzpicture}[main node/.style={inner sep=0,minimum size =.08cm,circle,fill=black!80,draw},node distance = 0.75cm,scale=0.5, baseline=0]
    \node[main node] (a) at (0,0) {};
    \node[main node] (d) at (2,0) {};
    \node[main node, label=above left:\small$x$] (g) at (-1,1) {};
    \node[main node] (h) at (-1,-1) {};
    \node[main node] (i) at (3,1) {};
    \node[main node] (j) at (3,-1) {};

    \draw (d) -- (a);
    \draw (g) -- (a) -- (h);
    \draw (i) -- (d) -- (j);
  \end{tikzpicture}
}

\newcommand{\ShapeOneCut}{%
  \begin{tikzpicture}[main node/.style={inner sep=0,minimum size =.08cm,circle,fill=black!80,draw},node distance = 0.75cm,scale=0.5, baseline=0]
    \node[main node] (a) at (0,0) {};
    \node[main node, label=above:\small$v$] (d) at (2,0) {};
    \node[style={inner sep=2,rectangle,draw}] (v) at (2,0) {};
    \node[main node, label=above left:\small$x$] (g) at (-1,1) {};
    \node[main node] (h) at (-1,-1) {};
    \node[main node] (i) at (3,1) {};
    \node[main node] (j) at (3,-1) {};

    \node (c1) at (1,1) {};
    \node (c2) at (1,-1) {};

    \draw[dashed] (c1) -- (c2);
    \draw (d) --  (a);
    \draw (g) -- (a) -- (h);
    \draw (i) -- (d) -- (j);
  \end{tikzpicture}
}

\newcommand{\ShapeTwoCut}{%
  \begin{tikzpicture}[main node/.style={inner sep=0,minimum size =.08cm,circle,fill=black!80,draw},node distance = 0.75cm,scale=0.5, baseline=0]
    \node[main node, label=above:\small$v$] (a) at (0,0) {};
    \node[style={inner sep=2,rectangle,draw}] (v) at (0,0) {};
    \node[main node] (d) at (2,0) {};
    \node[main node, label=above left:\small$x$] (g) at (-1,1) {};
    \node[main node] (h) at (-1,-1) {};
    \node[main node] (i) at (3,1) {};
    \node[main node] (j) at (3,-1) {};
    \node (c1) at (1,1) {};
    \node (c2) at (1,-1) {};

    \draw[dashed] (c1) -- (c2);

    \draw (d) -- (a);
    \draw (g) -- (a) -- (h);
    \draw (i) -- (d) -- (j);
  \end{tikzpicture}
}

\newcommand{\ShapeFive}{%
  \begin{tikzpicture}[main node/.style={inner sep=0,minimum size =.08cm,circle,fill=black!80,draw},node distance = 0.75cm,scale=0.5, baseline=0]
    \node[main node, label=above:\small$v$] (a) at (0,0) {};
    \node[main node] (g) at (-1,1) {};
    \node[main node] (h) at (-1,-1) {};
    \node[main node] (b) at (2,0) {};
    \node[main node] (c) at (2,-1.75) {};
    \node[main node, label=above:\small$w$] (d) at (4,0) {};
    \node[main node] (i) at (5,1) {};
    \node[main node] (j) at (5,-1) {};

    \draw (a) -- (b) -- (d);
    \draw (b) -- (c);
    \draw (g) -- (a) -- (h);
    \draw (i) -- (d) -- (j);
  \end{tikzpicture}
}

\newcommand{\ShapeFiveCut}{%
  \begin{tikzpicture}[main node/.style={inner sep=0,minimum size =.08cm,circle,fill=black!80,draw},node distance = 0.75cm,scale=0.5, baseline=0]
    \node[style={inner sep=2,rectangle,draw}] (v) at (0,0) {};
    \node[main node, label=above:\small$v$] (a) at (0,0) {};
    \node[main node] (g) at (-1,1) {};
    \node[main node] (h) at (-1,-1) {};
    \node[main node] (b) at (2,0) {};
    \node[main node, label=below:\small$x'$] (c) at (2,-1.75) {};
    \node[style={inner sep=2,rectangle,draw}] (w) at (4,0) {};
    \node[main node, label=above:\small$w$] (d) at (4,0) {};
    \node[main node] (i) at (5,1) {};
    \node[main node] (j) at (5,-1) {};
    \node (c1) at (1,1) {};
    \node (c2) at (1,-1) {};
    \node (c3) at (3,1) {};
    \node (c4) at (3,-1) {};

    \draw[dashed] (c1) -- (c2);
    \draw[dashed] (c3) -- (c4);

    \draw (a) -- (b) -- (d);
    \draw (b) -- (c);
    \draw (g) -- (a) -- (h);
    \draw (i) -- (d) -- (j);
  \end{tikzpicture}
}

%% file: intro.tex
\section{Introduction}\label{S:intro}

In 2009, in order to better understand Feynman integrals, Oliver Schnetz \cite{fq} defined an arithmetic graph invariant, the \textbf{$\mathbf{c_2}$-invariant}.  This invariant is not sensitive to all the information of the Feynman integral, and in particular many graphs are known which have the same $c_2$-invariant but which differ in their \textbf{Feynman period}, a particular important residue of the Feynman integral.  However, Brown and Schnetz conjectured \cite{modular}, based on now quite substantial computational evidence, that if two graphs have the same Feynman period then they have the same $c_2$-invariant.  If this is true, then known graph symmetries of the Feynman period must also be symmetries of the $c_2$-invariant.  One such symmetry the Feynman period is known to satisfy is the \textbf{completion symmetry}.  That is, if we begin with a 4-regular graph and remove one vertex, then the Feynman period does not depend on which vertex is removed.  In view of this, Brown and Schnetz conjectured in 2010 \cite{k3} that the $c_2$-invariant also satisfies completion symmetry.

This completion conjecture for the $c_2$-invariant has turned out to be surprisingly difficult to prove.  On the face of it the conjecture is a completely elementary statement involving only graph theory and elementary number theory, however, the motivic underpinnings connecting the Feynman integral and the $c_2$-invariant hint at its greater depth and difficulty.  Physicists and algebraic geometers both pursued the conjecture, but the first progress that did not depend on quite special graph structures such as double triangles, was the result of one of us \cite{specialc2} which proved the completion conjecture for $p=2$ when the initial 4-regular graph has an odd number of vertices.  This parity dependence was an unsatisfying artefact, and here we complete the proof for $p=2$ and make some comments on generalizing the approach to primes $p>2$.

While completing the writing of this paper, we learned of a new proof due to Erik Panzer of the completion conjecture for $c_2$ when $p=2$.  This proof uses completely different, but also combinatorial, tools and fundamentally cannot be generalized to $p>2$.  We look forward to the appearance of this new result.

The structure of the paper is as follows.
We begin Section~\ref{S:background} with our set up, including the definitions of the $c_2$-invariant and the completion conjecture, and more on the motivation for them.
We then define the polynomials we need, reduce the problem of computing $c_2$ to an enumerative problem, and lay out the three fundamental cases of our argument.  Sections~\ref{S:T-case},~\ref{S:S-case}, and \ref{S:R-case} then give the argument for $p=2$ in these three cases and the $p=2$ argument is completed in Section~\ref{S:completing-conjecture}.  We conclude the paper by presenting some partial results for $p>2$ in Section~\ref{S:higherp}.

\bpointn{Acknowledgements.}
KY is supported by an NSERC Discovery grant and by the Canada Research Chairs program.
The results of the present paper first appeared in the master's thesis of SH \cite{Hmmath}.\\

%% file: background.tex
\section{Background}\label{S:background}

Let $G$ be a connected graph and associate a Schwinger parameter $\alpha_e$ to each edge $e$ in $G$.  Then define the \textbf{graph polynomial} or \textbf{Kirchhoff polynomial} of $G$ to be
\[ \Psi_G = \sum\limits_{\substack{T \\ \text{spanning tree}}} \prod_{e \notin T} \alpha_e. \]
This polynomial is also known as the first Symanzik polynomial and is a particular variant of the spanning tree generating function.

In~\cite{fq}, Schnetz introduced the following arithmetic invariant and showed this invariant is well-defined by proving that $[\Psi_G]_q$ was indeed divisible by $q^2$ for $G$ connected with at least 3 vertices.

\tpointn{Definition} (Theorem $2.9$ of~\cite{fq})\label{def:c2}
\statement{
  Let $q = p^n$ be a prime power and $\mathbb{F}_q$ the finite field with $q$ elements.
  Let $G$ be a connected graph with at least $3$ vertices.
  Then the \textbf{$\mathbf{c_2}$-invariant} of $G$ at $q$ is
  \begin{equation}\label{eq:originalc2}
    c_2^{(q)}(G) \equiv \frac{\left[\Psi_G\right]_q}{q^2} \mod q,
  \end{equation}
  where $\left[\Psi_G\right]_q$ is the number of zeros of $\,\Psi_G$ in $\mathbb{F}_q^{\abs{E(G)}}$. \\
  Denote by $\mathbf{c_2(G)}$ the sequence of $c_2^{(q)}(G)$ for all prime powers $q$.
}

While the $c_2$-invariant is defined for all prime powers $q$ and more general connected graphs, we will restrict to when $q = p$ is a prime (see Conjecture~\ref{prime-powers}) and to decompletions of connected {4-regular} graphs $G$.
When specifically referring to the $c_2$ at primes $q=p$, we will often use $p$ instead of $q$.
When the $c_2$ is a constant sequence consisting of all $c$'s, we may say that $c_2 = c$.

The graphs which we are most interested in are those which are the result of taking a 4-regular graph and removing one vertex.  These graphs are most interesting to us because they can be interpreted as 4-point Feynman diagrams in scalar $\phi^4$ theory.  If the underlying 4-regular graph furthermore has no four edge cuts other than those which separate off a single vertex then we say the graph with the vertex removed is \textbf{primitive divergent}. The Feynman integrals of such graphs can be expressed in parametric form by using $\Psi_G$ and another polynomial.  Already very interesting is a particular residue of the integral called the \textbf{Feynman period}~\cite{periods} which is renormalization scheme independent and captures much of the number-theoretic content for massless scalar field theories in four-dimensional space-time.  The Feynman period in parametric form is defined as
\begin{equation}\label{eq:period}
  P_G \coloneqq \bigintsss_{\; 0}^{\infty} \frac{\;\mathrm{d}\alpha_1 \cdots \,\mathrm{d}\alpha_{\abs{E(G)}-1}\;}{\Psi_G^2 \,\bigg|_{\alpha_{\abs{E(G)}}=1}}.
\end{equation}
Note that this integral is independent of the choice of edge $e=\abs{E(G)}$ to set $\alpha_e = 1$.
There are many equivalent formulations of this integral, in the different spaces (position, momentum etc.) as well projective versions of each.
We refer the reader to ~\cite{census} and \textsection 3.1, 3.2 of~\cite{periods} for further details.

Looking at the period in its parametric form~\eqref{eq:period}, notice we are integrating over the denominator $\Psi_G^2$, which is just a polynomial in $\abs{E(G)}$ variables.
In particular, in order to understand and characterize properties of the period, we need to understand the structure of $\Psi_G$.
This motivates the study of the geometry of the variety defined by $\Psi_G$ and hence of the zeros of $\Psi_G$.

In the calculations of the 90s, see Broadhurst and Kreimer~\cite{knots,association}, all known Feynman periods were multiple zeta values, Kontsevich recognized that if the Feynman periods were all multiple zeta values because the underlying motives were mixed Tate, then the point count function $[\Psi_G]_q$ would be a polynomial in $q$ for each $G$, and so in 1997 he informally conjectured that this was the case.
It was since determined that not all Feynman periods are multiple zeta values and not all $[\Psi_G]_q$ are polynomial in $q$ \cite{BB, doryn, fq, k3}.
In cases when $[\Psi_G]_q$ is indeed a polynomial in $q$, $c_2^{(q)}(G)$ would then exactly be the quadratic coefficient of $[\Psi_G]_q$, hence the name "$c_2$".

In a sense, the $c_2$-invariant is measuring how badly Kontsevich's conjecture fails for a given graph.
More importantly, other than its combinatorial flavour which we will see shortly, the interest in the $c_2$-invariant lies in how it seems to detect the types of numbers that appear in the period and thus tells us something about the geometries underlying them~\cite{k3,modular,geometries}.  As mentioned above, Brown and Schnetz conjectured (Remark 2.11 (2) of~\cite{fq}, Conjecture $2$ of ~\cite{modular}) that if two primitive divergent graphs have the same period then they have the same $c_2$-invariant.
Currently, this conjecture holds for all known examples.

Finally, we note that we will restrict the $c_2$ sequence to all primes $q=p$ instead of prime powers $q=p^n$.
This is based on the following conjecture by Schnetz which says that prime powers do not give any extra information.

\tpointn{Conjecture} (Conjecture 2 of~\cite{geometries})\label{prime-powers}
\statement[eq]{
  Let $G_1$ and $G_2$ be two graphs with equivalent $c_2^{(p)}$ for all primes $p$, that is $c_2^{(p)}(G_1) \equiv c_2^{(p)}(G_2) \mod p$.
  Then for all prime powers $q = p^n$
  \[ c_2^{(q)}(G_1) \equiv c_2^{(q)}(G_2) \mod q. \]
}

In~\cite{census}, Schnetz proved that the period could be defined for the graphs created by adding back the vertex that was "removed" from primitive divergent graphs, and thus "completing" them to become 4-regular.
We call this 4-regular graph $G$ the \textbf{completion} and call the original primitive divergent graph $G-v$ for some vertex $v$ a \textbf{decompletion} of $G$.
Note that decompletions of 4-regular graphs could be non-isomorphic.
To utilize this, Schnetz further proved that if two primitive divergent graphs have the same completion, then their periods are also the same.
Thus the power of completion is that it allows us to lift the problem of calculating periods on primitive divergent graphs to one on 4-regular graphs, in the process cutting down the number of graphs to compute periods for.

For the $c_2$-invariant, whether or not the completion symmetry holds is still an open conjecture.
\tpointn[f]{The $c_2$ Completion Conjecture} (Conjecture $4,35$ of ~\cite{k3})\label{c2-comp-conj} \\
\statement[eq]{
  Let $G$ be a connected 4-regular graph, and let $v$ and $w$ be vertices of $G$.
  Then,
  \[ c_2(G - v) = c_2(G - w). \]
}

In~\cite{specialc2} one of us first made progress in the special case of $p=2$ by reducing the conjecture to a combinatorial counting problem; one involving enumerating certain edge bipartitions.
That paper partially proves the $p=2$ case using two different constructions.
The obstruction to a full proof lay in the second construction, which required restricting the argument to completed graphs with an odd number of vertices.
In this paper, inspired by the first construction in~\cite{specialc2}, we present this part of the argument and show how by using the same ideas we can get rid of the mysterious parity condition on the vertices, thereby completing the $c_2$ completion conjecture for this special case.

\bpoint{Graph polynomials}\label{S:graph-polys}

Before we can work on massaging Equation~\eqref{eq:originalc2}, we need to familiarize ourselves with the function in the numerator, the graph polynomial.
Notice that $\Psi_G$ is a homogeneous polynomial of degree $\ell(G) = |E(G)| - |V(G)| + 1$ in $|E(G)|$ variables.  We use $\ell(G)$ for the dimension of the cycle space of $G$, which is called the \textbf{loop number} in quantum field theory.
There is also a nice deletion-contraction relation, where for any edge $e$ in $G$ we have
\begin{equation}\label{eq:del-cont-psi}
  \Psi_G = \alpha_e \Psi_{G \setminus e} + \Psi_{G / e},
\end{equation}
since we can partition the sum based on whether a spanning tree of $G$ contains $e$ or not.

As the graphs get larger, these polynomials explode, and thus the first step in transforming the $c_2$ is to find smaller polynomials which we can take point counts of.
To do this, we need the theory of some related graph polynomials, which following Brown we call \textbf{Dodgson polynomials}, and which stem from the deletion-contraction relation of the graph polynomial and Kirchhoff's matrix-tree theorem.

The underlying framework of this theory on graph polynomials relies on the fact that the graph polynomial can be seen as the determinant of a particular matrix.
To represent the graph polynomial as such, we first define the following.

\tpointn{Definition}\label{def:L}
\statement{
  Given a connected graph $G$, choose an arbitrary orientation on the edges and an order $\iota$ on the edges and vertices of $G$, where the edges are ordered before the vertices,
  \[ \iota : E(G) \cup V(G) \quad\longrightarrow\quad \{1, \ldots, \abs{E(G)}+\abs{V(G)} \}.\]
  We denote by $\iota_e$ or $\ \iota_v$, to mean the order of edge $e$ or vertex $v$ under $\iota$, respectively.\\\\
  Let $\cE_G$ be the $\abs{V(G)} \times \abs{E(G)}$ signed incidence matrix, with any one row (corresponding to a vertex) removed.
  That is, for all vertices $v$ except one, and all edges $e$
  \[ [\cE_G]_{v,e} = \begin{cases} 1 & \text{ if }\, e = (v,u), u \in V(G), \\ -1 & \text{ if }\, e = (u,v), u \in V(G), \\ 0 & \text{ otherwise. }\end{cases} \]
  Let $\Lambda$ be the diagonal matrix of indeterminates $\alpha_e$ for $e$ in $E(G)$ in the chosen order, that is
  \[ [\Lambda]_{e_i,e_j} = \begin{cases} \alpha_{e_i} & \text{ if }\, i=j, \\ 0 & \text{ otherwise. }\end{cases} \]
  Then we define the \textbf{expanded Laplacian} of $G$ to be
  \[ L_G = \left[
    \renewcommand\arraystretch{1.5}
    \begin{array}{c|c}
      \Lambda & {\cE_G}^T \\
      \hline
      \cE_G & 0 \\
    \end{array}
    \right],
  \]
  which is an $(\abs{E(G)}+ \abs{V(G)} - 1) \times (\abs{E(G)} + \abs{V(G)} - 1)$ matrix, with rows and columns ordered by $\iota$.
}

While this matrix is not well-defined, as it depends on the choice of row removed in $\cE_G$ as well as the choice of orderings and orientation, we have, for any such choice
\begin{equation}\label{eq:det-psi}
  \Psi_G = (-1)^{|V(G)|-1}\det\left(L_G\right).
\end{equation}
This is a consequence of the matrix-tree theorem. See \textsection 2.2 of~\cite{periods} for a proof and \cite{geometries} for more discussion of the signs.

Under this determinantal framework and looking at the deletion-contraction relation~(\ref{eq:del-cont-psi}) for $\Psi_G$, we see that removing an edge $e$ to get $\Psi_{G\setminus e}$ corresponds to the minor where the row and column indexed by $e$ is removed.
This comes from noticing that $\Psi_{G\setminus e}$ is the coefficient of $\alpha_e$ in $\det(L_G)$.
Contracting an edge $e$ to get $\Psi_{G / e}$ then corresponds to setting $\alpha_e$ to zero.
In either case, some consideration is needed for the signs.
Motivated by this observation, we can extend the definition to minors of $L_G$ which we call \textbf{Dodgson polynomials}, or just Dodgsons.

\tpointn{Definition}\label{def:dodgson} (Definition 7 of~\cite{geometries})
\statement{
  Let $I$ and $J$ be words in the edges such that $\abs{I} = \abs{J}$ and let $K$ be a subset of edges of $G$.
  Denote $L_G{(I, J)}_K$ the matrix obtained from $L_G$ by removing rows indexed by $I$ and columns indexed by $J$, and setting $\alpha_e = 0$ for $e \in K$.
  Then the \textbf{Dodgson polynomial} is defined to be
  \[ \Psi_{G, K}^{I, J} = (-1)^{\abs{V(G)} + \iota_I + \iota_J -1} \sgn(I) \sgn(J) \det L_G(I, J)_K, \]
  We define $\Psi_{G,K}^{I,J} = 0$ if $\abs{I} \neq \abs{J}$ and define the empty determinant to be 1.  The specifics of the signs are not important for our purposes and so we leave their definitions to \cite{geometries}.
}

While there is a choice of ordering $\iota$ and vertex removed in $\cE_G$ for the matrix $L_G$, Dodgsons do not depend on this choice (see Lemma 9 of~\cite{geometries}).
When the graph $G$ is clear from the context, we will drop the subscript of $G$.
We also omit empty indices.
When $I, J, K$ are all empty, we recover $\Psi_G$.

Lastly, like how the monomials of $\Psi_G$ correspond to spanning trees of $G$, there is a nice combinatorial interpretation of the monomials appearing in Dodgsons.

\tpointn{Theorem}\label{dodgsonsum} (Proposition 23 of~\cite{periods})
\statement{
  Suppose $I \cap J = K = \emptyset$.
  Then we have
  \[ \Psi_{G}^{I, J} = (-1)^{\iota_I + \iota_J} \sgn(I) \sgn(J) \sum_{U \subset G \setminus (I \cup J)}  \det\left( \cE_G( G \setminus (U \cup I)) \right) \det\left( \cE_G( G \setminus (U \cup J)) \right)\prod_{e \not\in U} \alpha_e, \]
  where the sum runs over all subgraphs $U$ such that $U \cup I$ and $U \cup J$ are both spanning trees in $G$.
}

Ignoring the signs, what this theorem is saying is that we can view Dodgson polynomials through the possible shapes of the underlying graph after some deletions and contractions.
For $U \cup I$ to be a spanning tree in $G$, we must have that $U$ is a spanning tree in $G \setminus J / I$ since $U$ are edges not in $I$ and $J$.
Similarly for $U \cup J$ to be a spanning tree, $U$ must be a spanning tree in $G \setminus I / J$.
The polynomial can then be thought of as an "intersection" of graphs $\left(G \setminus J / I\right) \cap \left(G \setminus I / J\right)$ where the graphs represent the spanning trees and $\cap$ is taken to mean the resulting polynomial of common terms which are spanning trees in each minor.
When $K \neq \emptyset$, this corresponds to contracting the edges in $K$ (that are not in $I, J$) in both graphs.
When $I \cap J \neq \emptyset$, this corresponds to removing the edges in $I \cap J$ in both graphs.

However, here we are still dealing with the common spanning trees between two possibly different graphs.
In an effort to better understand these polynomials combinatorially, one of us with Brown introduced spanning forest polynomials~\cite{forest}.

\tpointn{Definition} (Definition 9 of~\cite{forest})\label{def:span-forest}
\statement{
  Let $P = P_1 \cup \cdots \cup P_k$ be a set partition of a subset of the vertices of $G$.
  Let $F$ be a spanning forest that partitions the vertices of $P$ exactly into $P_1 \cup \cdots \cup P_k$.
  More precisely, if $F = T_1 \cup \cdots \cup T_k$ then each tree $T_i$ of $F$ contains all the vertices in $P_i$ and no other vertices of $P$, and possibly other vertices in $V(G) \setminus P$.
  Then we say that the spanning forest $F$ is \textbf{compatible} with the vertex partition $P$.
  Note that the vertices not in $P$ can belong to any tree of $F$. \\\\
  We define a \textbf{spanning forest polynomial} of $G$ to be
  \[ \Phi_G^P = \sum_F \prod_{e \not\in F} \alpha_e, \]
  where the sum runs over all spanning forests $F$ that are compatible with $P$.
}

To relate spanning forest polynomials back to Dodgsons, instead of viewing Dodgsons through the matrix-tree theorem like in Theorem~\ref{dodgsonsum}, we can interpret them via the all-minors matrix-tree theorem~\cite{minorsmatrixtree}.
As a consequence, Dodgsons can be viewed as sums of spanning forest polynomials.

\tpointn{Theorem} (Proposition 12 of~\cite{forest})\label{tree-to-forest}
\statement{
  Let $I$, $J$ and $K$ be sets of edges of $G$ where $\abs{I} = \abs{J}$.
  Then
  \[ \Psi_{G,K}^{I,J} = \sum f_P \Phi_{G \setminus (I \cup J \cup K)}^{P}, \]
  where $f_P \in \{-1,1\}$, and the sum runs over all partitions $P$ of the vertices in $(I \cup J \cup K) \setminus (I \cap J)$ such that all forests compatible with $P$ are spanning trees in both $G \setminus I / (J \cup K)$ and $G \setminus J / (I \cup K)$.
}

\bpoint{Transforming the \texorpdfstring{$c_2$}{c2}-invariant}\label{S:denominator}

With the theory of graph polynomials in place, we return to our main goal of transforming the definition of the $c_2$-invariant~(\ref{eq:originalc2}) into a more combinatorial form.  We need two results.

Taking three edges $1,2,3$ of the graph $G$ in that order we can define
\[ D^3_G(1,2,3) = \pm \Psi_G^{13,23}\Psi^{1,2}_{G,3}. \]
Then,
\tpointn{Theorem}\label{denomred}
\statement[eq]{
  Let $G$ be a connected graph with $2\ell(G) \leq \abs{E(G)}$ and $\abs{E(G)} \geq 3$.
  Let $1,2,3$ be any three edges of $G$.
  Then
    \[ c_2^{(q)}(G) \equiv - \left[ \Psi_G^{13,23} \Psi_{G,3}^{1,2} \right]_q \mod q.\]
}

This theorem is a consequence of Corollary 28 of~\cite{k3} applied in the case $n=3$.  It is a special case of a much more general result that says that the denominators obtained by applying Brown's denominator reduction algorithm (see \cite{periods,massless}) to the Feynman period can each be used to calculate the $c_2$-invariant without any need to divide by $p^2$.  See Theorem $29$ of ~\cite{k3} with Corollary 28 of ~\cite{k3} for $n<5$.

While Theorem~\ref{denomred} brings us to a much more tractable form of the $c_2$-invariant, we are still left with the problem of point counting, which itself is a hard problem as the denominators can still have many variables.
To tackle this we turn to a consequence of a theorem from number theory, see \textsection 2 of~\cite{ax}, which will allow us to determine coefficients as a way of counting zeros.

\tpointn{Theorem} (Corollary of Chevalley-Warning Theorem)\label{chevalley}
\statement{
  Let $F$ be a polynomial of degree $N$ in $N$ variables, $x_1,\,\dots\, ,x_N$, with integer coefficients. Then we have
  \[ \text{coefficient of } x_1^{q-1}\,\cdots\, x_N^{q-1} \text{ in } F^{q-1} \equiv (-1)^{N-1}[F]_q \mod p, \]
  where $q$ is some power of the prime $p$.
}

Note that this theorem holds when reducing modulo prime $p$ and not when reducing modulo $q$.
An example of this is given in \textsection 5 of~\cite{geometries}.
Thus, to use this theorem in the context of the $c_2$-invariant, we will restrict to when $q=p$ is a prime (see Conjecture~\ref{prime-powers}).

To apply this theorem to $D^3_G$, we just need to make sure that the degree of the polynomial matches up with the number of variables, which it does since for a 4-regular graph with a vertex removed, $\Psi_G$ has degree exactly half the number of edges, then $\Psi_{G,3}^{1,2}$ has degree one less and $\Psi_G^{13,23}$ has degree two less, so $D^3_G$ has degree three less than the number of edges of $G$ and involves one variable for each edge of $G$ other than edges $1,2,3$.
Applying this to Theorem~\ref{denomred} and now using $[\cdot]$ to mean taking coefficients, as is common in enumerative combinatorics, we arrive at the second transformation for the $c_2$-invariant.

\tpointn{Corollary}
\statement[eq]{
  Let $G$ be a connected graph with $2\ell(G) \leq \abs{E(G)}$ and $\abs{E(G)} \geq 3$.
  Then, for any three distinct edges labelled $1,2,3$ of $G$ and when $\abs{E(G)} = 2\ell(G)$ we have
\begin{equation}\label{eq:c2-version}
  c_2^{(p)}(G) \equiv -[\alpha_4^{p-1}\cdots\alpha_{|E(G)|}^{p-1}]\; \left(\Psi_G^{13,23} \Psi_{G,3}^{1,2} \right)^{p-1} \mod p.
\end{equation}
}

What's nice about this version of the $c_2$-invariant is that we have reformulated what was a problem on counting zeros of polynomials into a problem about finding the coefficient of a particular monomial in some power of said polynomials.
Furthermore, these polynomials inherently have a combinatorial nature as they arise from spanning trees or spanning forests of particular graph minors.
Thus, we can interpret this coefficient as counting the number of ways to distribute edges into some of these spanning trees or spanning forests, which is a purely combinatorial problem!

Let us specialize further to the cases that we will need to attack the completion conjecture.
Let $v$ and $w$ be vertices of the connected 4-regular graph $G$.
Notice that it suffices to prove the completion conjecture holds when $v$ and $w$ are adjacent.
This is because as $G$ is connected, there is a path between any two vertices, and successively applying the equation to each pair of neighbours along the path gives the more general result.

\begin{figure}[t]
  \centering
  \includegraphics[scale=0.75]{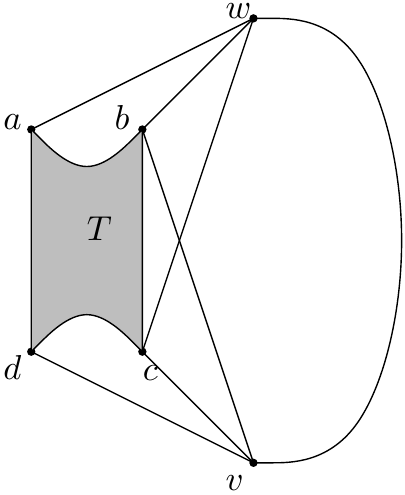}
  \quad\quad\quad \includegraphics[scale=0.75]{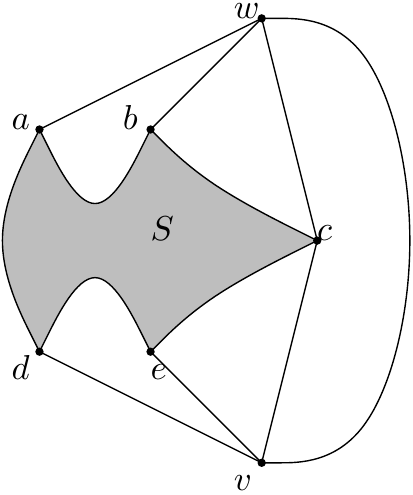}
  \quad\quad\quad \raisebox{4mm}[0pt][0pt]{\includegraphics[scale=0.75]{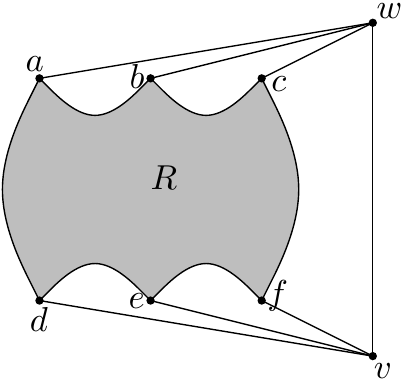}}
  \caption{The graph $G$, and the three cases $T$, $S$, and $R$.}
  \label{fig:TSR}
\end{figure}

As we will see shortly, the neighbours of $v$ and $w$ play an important role, and thus we split the conjecture into cases based on the number of neighbours that $v$ and $w$ have in common.
Since $G$ is a 4-regular graph and $v$ and $w$ are adjacent, there are four possibilities: they share all three neighbours, or exactly two or one neighbour, or they have no neighbours in common.
In the first case, when $v$ and $w$ share all neighbours, $G-v$ and $G-w$ are isomorphic, and thus the conjecture holds trivially.
For the others, we will deal with each case separately and label them with $T$, $S$, and $R$, respectively.
The $T$-case is when $v$ and $w$ share exactly two neighbours, the $S$-case is when they share only one neighbour, and the $R$-case is when they do not share any neighbours.
The three cases are depicted in Figure~\ref{fig:TSR}, where the grey blobs are the graphs $G - \{v,w\}$ in each case, and will be addressed in Sections~\ref{S:T-case},~\ref{S:S-case}, and~\ref{S:R-case}, respectively.

With the set-up complete, the first step in tackling the $c_2$ completion conjecture is the question of how to reduce it into an approachable combinatorial counting problem.

Now suppose $u$ is a 3-valent vertex in $G - v$.
Let the edges $1,2,3$ be the three edges incident to $u$, and let $u_1, u_2, u_3$ be the respective vertices incident to these edges as below:
\[ \Blob \]
Looking at the term $\Psi^{1,2}_3$, as discussed in Section~\ref{S:graph-polys} we can view this Dodgson polynomial as coming from the following intersection of sets of spanning trees
\[ \ExampleBlob \]
where the first blob depicts spanning trees on $(G -v) \setminus 1 / \{2,3\}$, and the second blob is on $(G - v) \setminus 2 / \{1,3\}$.

To translate this into spanning forests on $(G - v) \setminus \{1,2,3\}$ using Theorem~\ref{tree-to-forest}, we just need to determine the partitions $P$ of $\{u, u_1, u_2, u_3\}$ (as these are the vertices of edges $\{1,2,3\}$) such that all the spanning forests compatible with $P$ are spanning trees in the above intersection.
First, spanning trees on the first blob correspond to spanning forests such that $u$, $u_2$, and $u_3$ are in different trees, and $u_1$ can be in any of the three trees.
Similarily, spanning trees on the second blob correspond to spanning forests such that $u$, $u_1$, and $u_3$ are in different trees, and there is no restriction for which tree $u_2$ is in.
As these spanning forests have exactly three parts, the intersection of these blobs are spanning forests such that $u$ is in one tree, $u_3$ is in another tree, and $u_1$ and $u_2$ are in the last tree.
Furthermore, this is the only possible partition of $\{u, u_1, u_2, u_3\}$ that satisfies the required condition.

In particular, there is only one spanning forest polynomial in the sum, which is $\Phi^{\{u\},\{u_3\},\{u_1,u_2\}}_{(G - v) \setminus \{1,2,3\}}$.
However, notice that $u$ is an isolated vertex in $(G - v) \setminus \{1,2,3\}$ since edges $1,2,3$ are exactly the edges incident to $u$.
This means we can remove the vertex $u$, giving the equation
\[ \Psi^{1,2}_{3} = \Phi^{\{u_3\},\{u_1,u_2\}}_{G - \{u,v\}}. \]
Thus we are left with spanning forests with exactly two parts, which we call \textbf{spanning 2-forests}, such that they are compatible with the vertex partition $\{u_3\},\{u_1,u_2\}$.

Following the same process for $\Psi^{13,23}$, in either case we want spanning trees on the graph $G - \{u,v\}$
\[ \BlobBase \]
where for $(G - v) \setminus \{1,3\} / 2$, we have $u = u_2$, and for $(G - v) \setminus \{2,3\} / 1$, we have $u = u_1$.
Once again, we can translate these to spanning forests on $(G - v) \setminus \{1,2,3\}$, this time with exactly two parts, giving the polynomial $\Phi^{\{u\},\{u_1,u_2,u_3\}}_{(G - v) \setminus \{1,2,3\}}$.
Noticing that $u$ is an isolated vertex in this graph, removing vertex $u$ gives
\[ \Psi^{13,23} = \Phi^{\{u_1,u_2,u_3\}}_{G - \{u,v\}} = \Psi_{G - \{u,v\}}. \]

Putting everything together, we have reformulated the $c_2$-invariant to \\
\begin{equation}\label{eq:c2-higher}
  c_2^{(p)}(G-v) \equiv -[\alpha_4^{p-1}\cdots\alpha_{|E|}^{p-1}]\; \left(\Phi_{G-\{u,v\}}^{\{u_3\},\{u_1,u_2\}} \Psi_{G-\{u,v\}} \right)^{p-1} \mod p,
\end{equation}

which we call \textbf{reducing with respect to vertex $u$}.
Note that the ordering of the vertices $u_1$, $u_2$, and $u_3$ was completely arbitrary, and this formulation of the $c_2$ does not depend on this choice.

Finally, by working with the $c_2$ in this form, there's a particularly nice interpretation of what this equation is counting.
From the coefficient extraction, the monomial we are looking for can be viewed as distributing $p-1$ copies of each edge in $G - \{u,v\}$ across $2p-2$ polynomials, where now $G - \{u,v\}$ is the underlying graph for the polynomials.
Of these $2p-2$  polynomials, $p-1$ of them arise from spanning 2-forests compatible with the vertex partition $\{u_3\},\{u_1,u_2\}$ and the other $p-1$ come from spanning trees.
Then, since all the monomials in $\Psi$ and $\Phi$ appear with a coefficient of $1$, we are enumerating all the ways to do such a distribution.

Now, recalling that the terms in $\Psi$ and $\Phi$ correspond to the edges \emph{not} in the respective trees or forests, we need to shift our view slightly to make sense of what this coefficient is counting.
Given $p-1$ spanning trees of $G - \{u,v\}$, the corresponding monomial in $\left(\Psi_{G - \{u,v\}}\right)^{p-1}$ comes from combining all the edges of the graph not in each spanning tree.
Thus the edges left over to be distributed to the $p-1$ $\Phi$ polynomials are exactly the edges of the $p-1$ spanning trees.
Furthermore, these edges must be grouped in a such a way that removing them gives $p-1$ spanning 2-forests.

Similarily, given $p-1$ spanning 2-forests the corresponding monomial in $\left(\Phi_{G-\{u,v\}}^{\{u_3\},\{u_1,u_2\}}\right)^{p-1}$ combines all the edges of the graph not in each forest.
Notice this monomial must then match the edges left over from the monomial given by the $p-1$ spanning trees above, which are exactly the edges of the spanning trees.
By the same logic, the monomial arising from the $p-1$ spanning trees must then match the edges of the $p-1$ spanning 2-forests.
Thus, we can swap the roles of the spanning trees and forests, and equivalently express the enumeration as counting the number of ways to partition $p-1$ copies of each edge into $p-1$ spanning trees and $p-1$ spanning 2-forests compatible with $\{u_3\},\{u_1,u_2\}$.

When $p=2$, this interpretation is much simpler as the equation simplifies to \\
\begin{equation}\label{eq:c2-counting}
  c_2^{(2)}(G-v) \equiv [\alpha_4\cdots\alpha_{|E|}]\; \Phi_{G-\{u,v\}}^{\{u_3\},\{u_1,u_2\}} \Psi_{G-\{u,v\}} \mod 2.
\end{equation}

Now we only have one copy of each edge, and we are counting the number of edge bipartitions where one part is a spanning tree and the other is a spanning 2-forest compatible with $\{u_3\},\{u_1,u_2\}$.
We will usually denote these edge partitions as $(\psi, \phi)$.
Additionally, as we are counting modulo $2$, we really only care about the parity of this count.
Thus one strategy to tackling the $c_2$ completion conjecture when $p=2$ is to find fixed-point free involutions on the appropriate sets of edge bipartitions! \\

%% file: t-case.tex
\section{The \texorpdfstring{$T$}{T}-case}\label{S:T-case}

Let $G$ be a connected 4-regular graph, and let $v$ and $w$ be two adjacent vertices of $G$ such that they share two common neighbours.
Let $T = G - \{v,w\}$ be the graph obtained by removing vertices $v$ and $w$ from $G$, the first grey blob in Figure~\ref{fig:TSR}, with the neighbours of $v$ and $w$ labelled $\{a,b,c,d\}$ as in the figure.
Here $w$ has neighbours $\{a,b,c,v\}$, and $v$ has neighbours $\{b,c,d,w\}$.

What's special about this case is when $v$ and $w$ share two common neighbours, these four vertices form a double triangle, which we know we can reduce as the $c_2$ is invariant under a transformation called \textbf{double triangle reduction} (see~\cite{k3,forest}).
Now adding in decompleting at a double triangle vertex, as proved in ~\cite{further}, we can reduce the $c_2$ completion conjecture in the $T$-case for all values of $p$ to the $c_2$ completion conjecture in the other cases on smaller graphs as a consequence of double triangle reduction.

This result can alternately be proved in an enumerative way, analogous to the arguments for the $S$ and $R$ cases below.  See \cite{Hmmath} for details.  We will return to this idea as a model for how an argument for $p>2$ could look in Section~\ref{S:higherp}.  Through whichever technique we obtain the following result.

\tpointn{Corollary}\label{T-case}
\statement[eq]{
  Let $G$ be a connected 4-regular graph.
  Let $v$ and $w$ be adjacent vertices of $G$ such that they share two common neighbours.
  Then,
  \[ c_2^{(2)}(G - v) = c_2^{(2)}(G - w). \]
} 

%% file: s-case.tex
\section{The \texorpdfstring{$S$}{S}-case}\label{S:S-case}

Let $G$ be a connected 4-regular graph, and let $v$ and $w$ be two adjacent vertices of $G$ such that they share one neighbour.
Let $S = G - \{v,w\}$ be the graph obtained by removing vertices $v$ and $w$ from $G$, the middle grey blob in Figure~\ref{fig:TSR}, with the neighbours of $v$ and $w$ labelled $\{a,b,c,d,e\}$ as in the figure.
Here $w$ has neighbours $\{a,b,c,v\}$, and $v$ has neighbours $\{c,d,e,w\}$.

As $v$ and $w$ still share a neighbour, after removing both vertices $S$ will have a 2-valent vertex.  For part of the argument we will be able to use the idea of swapping edges around the 2-valent vertex. This idea can then be generalized to swapping around a particular vertex called the control vertex, under certain conditions.
While first used in ~\cite{specialc2} for the $R$-case, we prove that this control vertex argument can be modified to apply in the $S$-case as well and thus complete the proof of the $c_2$ completion conjecture for $p=2$ in the $S$-case.

\bpoint{Set-up}

We start with defining the particular sets of edge partitions that we are counting.

\tpointn{Definition} (Definition 3.1 of~\cite{specialc2})\label{S-defs}
\statement{
  Suppose $P$ is a bipartition of $\{a,b,c,d,e\}$. \\
  Let $\mathcal{S}_P$ be the set of bipartitions $(\psi, \phi)$ of the edges of $S$ such that $\psi$ is a spanning tree and $\phi$ is a spanning 2-forest compatible with $P$.
  Let $s_P = \abs{\mathcal{S}_P}$.
}

We can pick the vertex bipartitions for the spanning 2-forests to exploit the property that $v$ and $w$ share a neighbour.

\tpointn{Proposition} (Proposition 3.2, 3.3 of~\cite{specialc2})\label{S-eqs}
\statement[eq]{
  When $v$ and $w$ have one common neighbour $c$
  \[ c_2^{(2)}(G - v) = s_{\{c\},\{a,b\}} \mod 2, \]
  \[ c_2^{(2)}(G - w) = s_{\{c\},\{d,e\}} \mod 2, \]
  and thus we have
  \begin{align*}\refstepcounter{equation}\tag{\theequation}\label{eq:S-counts}
    c_2^{(2)}(G - v) - c_2^{(2)}(G - w)
      =\ &{s}_{\{c,d\},\{a,b,e\}} + {s}_{\{c,e\},\{a,b,d\}} + {s}_{\{a,b\},\{c,d,e\}} \\
      &+ {s}_{\{a,c\},\{b,d,e\}} + {s}_{\{b,c\},\{a,d,e\}} + {s}_{\{d,e\},\{a,b,c\}} \mod 2.
  \end{align*}
}

The first part is the argument from Section~\ref{S:denominator}.
The second part follows from the first by summing over all possibilities for assigning the remaining vertices of $\{a,b,c,d,e\}$ to the parts of the partitions and then cancelling modulo 2.

\bpoint{Swapping around \texorpdfstring{$c$}{c}}

The main property of the $S$-case is that $v$ and $w$ share a neighbour, and in particular, after removing both $v$ and $w$ from $G$, we are left with a 2-valent vertex.  In all of the partitions of interest $c$ is in a part with at least one other vertex and so must be connected by at least one edge to the rest of the graph in both the spanning tree and the 2-forest.  However since $c$ is 2-valent, $c$ must be a leaf in both the spanning tree and the 2-forest.  Swapping which edge incident to $c$ is in the spanning tree and which is in the 2-forest, then, cannot create any cycles in either; the only change it can create is to switch which part of the partition contains $c$.  Using this idea we can prove the following theorem.

\tpointn{Theorem} (Lemma 4.2 of~\cite{specialc2})\label{S-swapc}
\statement[eq]{
  There is a fixed-point free involution on
  \[ \mathcal{S}_{\{a,b,c\},\{d,e\}} \cup \mathcal{S}_{\{a,b\},\{c,d,e\}}, \]
  and thus we have
  \[ {s}_{\{a,b,c\},\{d,e\}} + {s}_{\{a,b\},\{c,d,e\}} \equiv 0 \mod 2. \]
}

Swapping around $c$ takes care of two specific vertex bipartitions on the right-hand side of Equation~\eqref{eq:S-counts}, and what we are left with are sets of the form $\mathcal{S}_{\{c,*\}, \{*,*,*\}}$.
However, when we try to swap around $c$ for these sets we obtain vertex bipartitions not in the equation nor ones that coincide with each other.

\bpoint{Swapping around a control vertex}\label{SS:S-new}

For the remaining sets in Equation~\eqref{eq:S-counts}, we generalize the swapping argument to one using a vertex that is not 2-valent in $S$.
While this idea of using control vertices is not new, we adapt it to be able to apply it to the $S$-case.
In the process we develop a new swapping method involving multiple vertices, starting with using a control vertex in conjunction with a 2-valent vertex.
Generalizing once again, the notion of using multiple control vertices will then be used for the $R$-case.

\begin{figure}[t]
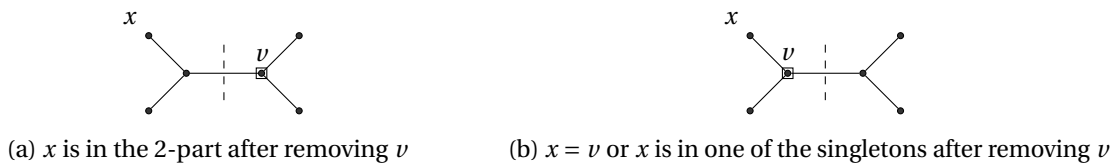

  \centering
  \begin{subfigure}[b]{0.5\linewidth}
    \centering
    \ShapeOneCut%
    \caption{$x$ is in the 2-part after removing $v$}
    \label{fig:S-control}
  \end{subfigure}%
  \begin{subfigure}[b]{0.5\linewidth}
    \centering
    \ShapeTwoCut%
    \caption{$x=v$ or $x$ is in one of the singletons after removing $v$}
    \label{fig:R-control}
  \end{subfigure}%
  \caption[The two possible control vertices for Lemma~\ref{control-vertex}.]{The two possible control vertices $v$ depending on the specification for special vertex $x$ for Lemma~\ref{control-vertex}. The three unlabelled leaves are the vertices in $p \setminus \{x\}$. The dashed lines indicate the edges we will be picking for one of the involutions in the $S$-case and in the $R$-case, respectively. Some of the leaf edges are possibly contracted, with the exception that two leaf edges incident to the same vertex cannot simultaneously both be contracted.}
  \label{fig:control}
\end{figure}

\tpointn{Lemma} (generalizing Lemma 4.3 of~\cite{specialc2})\label{control-vertex}
\statement{
  Let $G$ be a connected 4-regular graph with two adjacent vertices removed, and let $V$ be the set of neighbours of those vertices.
  Fix a vertex $x \in V$, which we call the \textbf{special vertex}.
  Suppose we have a bipartition $P$ of $V$ such that  one part is of the form $p = \{x, *, *, *\}$. \\\\
  Let $(\psi, \phi)$ be an edge partition of $G$ where $\psi$ is a spanning tree and $\phi$ is a spanning 2-forest compatible with $P$.
  Let $t$ be the tree corresponding to $p$ in the 2-forest $\phi$.
  Then, there are two distinct vertices $v$ such that removing $v$ from $t$ partitions $\{x,*,*,*\}$ into a 2-part, a singleton, and possibly a third part which is also a singleton. \\\\
  If we further specify that $x$ must satisfy one of the following properties, either
  \begin{enumerate}
    \item $x$ is in the 2-part, or
    \item $x=v$ or $x$ is in one of the singleton components,
  \end{enumerate}
  then there is a unique vertex $v$, which we call the \textbf{control vertex}, that satisfies all the required properties.
  Figure~\ref{fig:control} depicts the control vertex in each case.
  Furthermore, in $t$ this control vertex $v$ will either be 2-valent with $v \in p$ or 3-valent and $v \not\in p$.
  In both cases $v$ is a leaf in the spanning tree $\psi$.
}
\begin{proof}
  Looking at the tree $t$ and the subtree created by taking the union of all the paths in $t$ between any two vertices in $p$, since $p$ has four vertices this subtree will be of the form
  \[ \ShapeOne \]\\
  in which the leaves are vertices in $p$, the edges are paths in $t$, and some paths are possibly contracted.

  First, we argue that this is indeed the correct form for the subtree.
  Since $\psi$ is a spanning tree and we are partitioning the edges of $G$, the vertices in $\phi$, and thus the non-leaves in the subtree, are at most 3-valent.
  Furthermore, in $\phi$ the vertices in $P$, and thus in $p$, are at most 2-valent.
  This immediately implies that we have the subtree form above, where some of the paths may be contracted.
  The valency constraints also mean that the path between the two non-leaf vertices cannot be contracted.
  Note trivially, the paths between two leaves adjacent to the same closest non-leaf vertex cannot both be contracted simultaneously.
  All possible configurations are drawn in Figure 6 of~\cite{specialc2}.

  Now, because of the restrictions on the subtree configurations, removing any one of the two non-leaf vertices in the subtree above from $t$ creates components that partitions $p$ into a 2-part, a singleton, and possibly a third singleton.
  By the valency of the vertices in $\phi$, these two vertices can only be 2-valent if they are in $p$ or otherwise 3-valent.
  In either case, by the properties of $G$, they will be leaves in $\psi$.

  Let $v$ be one of the two non-leaf vertices in the subtree above.
  Finally, to get uniqueness of $v$ we need to fix where one vertex $x \in p$ must be after removing $v$ from $t$.
  In the first case, fixing $x$ to be in the 2-part component, we get the control vertex $v$ as pictured in Figure~\ref{fig:S-control}.
  Notice here $x$ could never be the control vertex $v$.
  In the second case, if we want $x$ to be in a singleton component we choose the other non-leaf vertex to be $v$ as in Figure~\ref{fig:R-control}.
  However, with this choice of $v$, we could also have that $x=v$ (the path between $x$ and $v$ is contracted in $t$).
  Thus in this case, the property we are fixing is either $x=v$ or $x$ is in one of the singleton components.
\end{proof}

Now we are ready to tackle the remaining sets in the Equation~\eqref{eq:S-counts}, which all are of the form $\mathcal{S}_{\{c,*\},\{*,*,*\}}$.
Instead of dealing with these sets directly, since we know there is a fixed-point free involution via swapping around $c$ (as it is a 2-valent vertex), we can equivalently show that there is a fixed-point free inovlution on sets of the form $\mathcal{S}_{\{*\},\{c,*,*,*\}}$ (the "swapped" version).

The involution uses a two-phase process, where we utilize that we can swap around $c$ for certain edge bipartitions.
For the others, we modify the involution used in Theorem~\ref{R-control-bij}, which uses the same idea of swapping edges incident to the control vertex between the two parts of $(\psi, \phi)$.

\tpointn{Theorem}\label{S-bijection}
\statement[eq]{
  There is a fixed-point free involution on
  \[ \mathcal{S}_{\{a\},\{b,c,d,e\}} \cup \mathcal{S}_{\{b\},\{a,c,d,e\}} \cup \mathcal{S}_{\{d\},\{a,b,c,e\}} \cup \mathcal{S}_{\{e\},\{a,b,c,d\}}, \]
  and thus we have
  \[ {s}_{\{a\},\{b,c,d,e\}} + {s}_{\{b\},\{a,c,d,e\}} + {s}_{\{d\},\{a,b,c,e\}} + {s}_{\{e\},\{a,b,c,d\}} \equiv 0 \mod 2. \]
}
\begin{proof}
  Consider the edge partitions in the union of sets of the form $\mathcal{S}_{\{*\},\{c,*,*,*\}}$ where the $*$'s are the vertices $\{a,b,d,e\}$ in any order.
  Since partitioning $\{a,b,c,d,e\}$ (the neighbours of the two adjacent vertices removed to obtain the graph $S$) into the form $\{*\},\{c,*,*,*\}$ leaves each part containing at least one of $\{a,b,c,d,e\}$ and a 4-part containing $c$, these edge partitions satisfy the conditions allowing us to swap edges around the 2-valent vertex (as described before Theorem~\ref{S-swapc} and further detailed below) and the conditions of
  Lemma~\ref{control-vertex}.
  For Lemma~\ref{control-vertex}, we are taking $c$ to be the special vertex and specifying that $c$ be in the 2-part after removing the control vertex from the appropriate tree in the spanning 2-forest (the first property).

  Let $(\psi, \phi)$ be an edge partition in any set of the form $\mathcal{S}_{\{x\},\{c,*,*,*\}}$, and let $v$ be its control vertex.
  Let $t$ be the tree in $\phi$ corresponding to part $\{c,*,*,*\}$ which also contains $v$.
  Since $v$ is a leaf in the spanning tree $\psi$, there is exactly one edge incident to $v$ in $\psi$, which we will call $\eta_v$, and let $n_v$ be the neighbour of $v$ in $\psi$.
  Since $c$ is a 2-valent vertex, and thus also a leaf in $\psi$, let $\eta_c$ be the edge incident to $c$ in $\psi$ and $n_c$ its neighbour.

  Now to describe the involution giving a new edge partition:
  \begin{enumerate}[label=(\arabic*)]
    \item\label{S-step-1} \emph{Swapping $c$ stays in} -- If $n_c \in t$, then swap the edges incident to $c$ between $\psi$ and $\phi$.
    \item\label{S-step-2} \emph{Swapping $c$ goes out} -- Otherwise:
      \begin{enumerate}[label=(\roman*), ref=(2\roman*)]
        \item\label{S-step-in} \emph{Swapping $v$ stays in} -- If $n_v \in t$, so the control vertex and its neighbour in $\psi$ are in the same tree of $\phi$, let $\eta$ be the edge incident to $v$ in $t$ in the path to $n_v$.\\
          Then, swap the edges $\eta_v$ and $\eta$ between $\psi$ and $\phi$.
        \item\label{S-step-out} \emph{Swapping $v$ goes out} -- If $n_v \not\in t$, so the control vertex and its neighbour in $\psi$ are in different trees of $\phi$, let $\eta$ be the edge incident to $v$ in $t$ in the path to $c$; in Figure~\ref{fig:S-control} this is the edge towards the 2-part indicated by the dashed line. \\
          Then swap the edges $\eta_v$ and $\eta$ between $\psi$ and $\phi$, and finally swap the edges incident to $c$.
      \end{enumerate}
  \end{enumerate}

  First, we show that we get a valid edge partition in a set of the form $\mathcal{S}_{\{*\},\{c,*,*,*\}}$.
  For ~\ref{S-step-1}, this is just a specific case of the swapping around a 2-valent vertex argument.
  Since $c$ is a 2-valent vertex and is a leaf in both $\psi$ and $\phi$, we can swap the edges incident to $c$ between the two parts to get $(\psi', \phi')$ where $\psi'$ is still a spanning tree and $\phi'$ is still a spanning 2-forest.
  Looking at $\phi'$, as $n_c \in t$ and thus both neighbours of $c$ were in the same tree in $\phi$, $c$ gets reconnected to $t$ after the swap.
  Thus $\phi'$ is compatible with the same vertex bipartition as $\phi$, and $(\psi', \phi')$ is in the same set $\mathcal{S}_{\{x\},\{c,*,*,*\}}$ as before.

  For ~\ref{S-step-2}, we can use a similar argument to swapping around a 2-valent vertex, except replacing the 2-valent vertex with the control vertex $v$.
  From Lemma~\ref{control-vertex}, we know that $v$ is the unique vertex such that removing $v$ from $t$ partitions the vertices in part $\{c,*,*,*\}$ into a 2-part containing $c$, a singleton, and possibly a third singleton.
  Unlike the 2-valent argument, as $v$ is no longer a leaf in $t$ and we now have a choice of which edge incident to $v$ in $t$ we are swapping, we need to be careful of creating cycles when reconnecting $v$ to the rest of $t$ and make sure that we can recover the edge swapped and the same control vertex to indeed get an involution.

  To deal with the possibility of creating cycles, we notice that the only way to get a cycle after swapping edges incident to $v$ between $\psi$ and $\phi$ is if the neighbour of $v$ in $\psi$ was in the same tree as $v$ in $\phi$ and the edge picked to be swapped in $\phi$ was not on the cycle created by adding the edge incident to $v$ in $\psi$ to $\phi$.
  Thus for ~\ref{S-step-in}, if $n_v \in t$, then we pick the edge $\eta$ to be the edge incident to $v$ on the path to $n_v$ in $t$.
  In this case all the neighbours of $v$ are in the same tree in $\phi$.
  Removing $\eta$ from $t$ breaks the tree $t$ into two components, one with $v$ and one containing $n_v$, splitting $\phi$ into three components.
  Adding $\eta_v$ reconnects the two components from $t$ via $v$ and $n_v$ to create a new spanning 2-forest $\phi'$.
  As $n_v$ and $v$ were both in $t$, $\phi'$ is compatible with the same vertex bipartition as $\phi$.
  Since $v$ is a leaf in $\psi$, removing $\eta_v$ from $\psi$ and adding $\eta$ to create $\psi'$ maintains the spanning tree structure.
  Thus swapping the edges $\eta_v$ and $\eta$ between $\psi$ and $\phi$ creates a new edge partition $(\psi', \phi')$ in the same set $\mathcal{S}_{\{x\},\{c,*,*,*\}}$ as we started with.

  For ~\ref{S-step-out} when $n_v \not\in t$, we pick $\eta$ to be the edge incident to $v$ on the path to $c$ in $t$.
  Removing $\eta$ from $t$ splits the tree into a component with $v$ and a component which contains $c$ and another vertex in part $\{c,*,*,*\}$, call it $y$, breaking $\phi$ into three components.
  Since $n_v \not\in t$, adding $\eta_v$ connects the component with $v$ to the other tree in $\phi$ (corresponding to part $\{x\}$) creating a new spanning 2-forest $\psi'$.
  Once again as $v$ is a leaf in $\psi$, removing $\eta_v$ from $\psi$ and adding $\eta$ to create $\psi'$ maintains the spanning tree structure.
  Thus swapping the edges $\eta_v$ and $\eta$ between $\psi$ and $\phi$ creates a new edge partition $(\psi', \phi')$,
  However, now $\psi'$ is compatible with the vertex bipartition $\{c,y\},\{x,*,*\}$.
  To get back to a vertex bipartition of the form $\{*\},\{c,*,*,*\}$, we notice that we are in the case where $n_c \not\in t$, and in particular, $n_c$ was in the other tree of $\phi$ which is now connected to $x$.
  Since $c$ can never be the control vertex, swapping the edges incident to $c$ between $\psi'$ and $\phi'$ gives a new edge partition $(\psi'', \phi'')$ now in the set $\mathcal{S}_{\{y\},\{c,x,*,*\}}$, which is of the correct form.
  Note that $y \neq x$ could be any of $\{a,b,d,e\}$.

  In all three cases, as we are always changing the edges incident to $c$ and/or $v$ in $\psi$ and $\phi$, the new edge partitions can never be identical to $(\psi, \phi)$.
  Thus this transformation is fixed-point free.

  Lastly, we prove that this transformation is indeed an involution by showing that the control vertex and the two edges being swapped remains the same in the new edge partition.
  Let $(\psi', \phi')$ be the transformed edge partition.
  For~\ref{S-step-1}, since $c$ is a 2-valent vertex, there is no choice of which edge is to be swapped and as both $n_c$ and $c$ are in the same tree in $\phi'$, we are in the same case of the transformation.
  Thus applying the transformation again returns us to $(\psi, \phi)$.

  For~\ref{S-step-2}, it will be easier to first formulate the control vertex in another way.
  We can define the control vertex $v$ to be the first common vertex in $t$ in the paths from exactly two of the vertices in part $\{c, *, *, *\} \setminus \{c\}$ to $c$.
  Under this formulation, we immediately see that any reconnecting of $v$ via $\eta_v$ does not change this property, and thus the control vertex remains the same after the swapping of edges incident to $v$.
  Additionally for ~\ref{S-step-out}, the swapping of edges incident to $c$ afterwards also does not affect this property of $v$ since $n_c$ was in the tree connected to $x$, and thus once again the control vertex remains the same in $\psi'$.

  For the edges, notice that~\ref{S-step-in} occurs when all neighbours of $v$ are in the same tree of $\phi$, and thus by the choice of $\eta$, all the neighbours of $v$ are still in the same tree of $\phi'$.
  Additionally in $\psi'$, $\eta_v$ is exactly the edge incident to $v$ on the cycle created by adding $\eta$, which is now the edge incident to $v$ in $\psi'$.
  For~\ref{S-step-out}, since we are reconnecting $v$ to the other tree of $\phi$ via $\eta_v$, the new neighbour of $v$ in $\psi'$ is now in a different tree than $v$ in $\phi'$.
  Additionally in $\psi'$, because $n_c$ was in the tree connected to $x$ in $\phi$, after the two swaps $\eta$ is exactly the edge incident to $v$ on the path to $c$.
  In both cases $(\psi', \phi')$ remain in the same case of ~\ref{S-step-2} as before with $\eta'_v = \eta$, as $v$ is a leaf in $\psi'$ and so there is only one choice of edge, and $\eta' = \eta_v$ as shown above.
  Thus, applying the transformation again returns us to $(\psi, \phi)$.

  Putting everything together, we indeed have a fixed-point free involution in all cases.
  Finally, this holds for any edge partition in sets of the form $\mathcal{S}_{\{x\},\{c,*,*,*\}}$ where we are either staying in the same set $\mathcal{S}_{\{x\},\{c,*,*,*\}}$ or swappping to a different set $\mathcal{S}_{\{y\},\{c,*,*,*\}}$ where $y \neq x$ could be any of $\{a,b,d,e\}$.
  Thus our involution is on the union of sets of the form $\mathcal{S}_{\{*\},\{c,*,*,*\}}$,
  \[ \mathcal{S}_{\{a\},\{b,c,d,e\}} \cup \mathcal{S}_{\{b\},\{a,c,d,e\}} \cup \mathcal{S}_{\{d\},\{a,b,c,e\}} \cup \mathcal{S}_{\{e\},\{a,b,c,d\}}, \]
  and therefore, the size of this set must be even, giving the equation
  \[ {s}_{\{a\},\{b,c,d,e\}} + {s}_{\{b\},\{a,c,d,e\}} + {s}_{\{d\},\{a,b,c,e\}} + {s}_{\{e\},\{a,b,c,d\}} \equiv 0 \mod 2. \]
\end{proof}

As in our discussion before Theorem~\ref{S-bijection}, since swapping around $c$ gives fixed-point free involutions $\mathcal{S}_{\{x\},\{c,*,*,*\}} \cup \mathcal{S}_{\{c,x\},\{*,*,*\}}$ for $x \in \{a,b,d,e\}$, we immediately get
\[ {s}_{\{a,c\},\{b,d,e\}} + {s}_{\{b,c\},\{a,d,e\}} + {s}_{\{c,d\},\{a,b,e\}} + {s}_{\{c,e\},\{a,b,d\}} \equiv 0 \mod 2 \]
by adding all the equations together and using the result from Theorem~\ref{S-bijection}.

Alternatively, notice that we could use the same fixed-point free involution on the union of sets of the form $\mathcal{S}_{\{c,*\},\{*,*,*\}}$ except with a slight modification.
In~\ref{S-step-2} where the two neighbours of $c$ are in different trees in the spanning 2-forest, we first swap around the edges of $c$, giving a new spanning 2-forest compatible with $\{*\},\{c,*,*,*\}$, and proceed as before with finding the control vertex and edges to swap.
Instead of swapping around $c$ again in ~\ref{S-step-out}, we do so in ~\ref{S-step-in} to ensure we get edge partitions in sets of the form $\mathcal{S}_{\{c,*\},\{*,*,*\}}$.
As everything else remains the same, we obtain the following result.

\tpointn{Corollary}\label{S-bij-cor}
\statement[eq]{
  There is a fixed-point free involution on
  \[ \mathcal{S}_{\{a,c\},\{b,d,e\}} \cup \mathcal{S}_{\{b,c\},\{a,d,e\}} \cup \mathcal{S}_{\{c,d\},\{a,b,e\}} \cup \mathcal{S}_{\{c,e\},\{a,b,d\}}, \]
  and thus we have
  \[ {s}_{\{a,c\},\{b,d,e\}} + {s}_{\{b,c\},\{a,d,e\}} + {s}_{\{c,d\},\{a,b,e\}} + {s}_{\{c,e\},\{a,b,d\}} \equiv 0 \mod 2. \]
}

\bpoint{Completing the \texorpdfstring{$S$}{S}-case}

With the two swapping arguments, first with swapping around $c$ and then with swapping around a control vertex, we have covered all the $\mathcal{S}$ sets in Equation~\eqref{eq:S-counts}!
In fact they are really the same argument as we can think of the 2-valent vertex $c$ as acting as a special type of control vertex.
Thus we have everything we need to prove the $c_2$ completion conjecture for $p=2$ in the $S$-case.

\tpointn{Theorem}\label{S-case}
\statement[eq]{
  Let $G$ be a connected 4-regular graph.
  Let $v$ and $w$ be adjacent vertices of $G$ such that they share one common neighbour.
  Then,
  \[ c_2^{(2)}(G - v) = c_2^{(2)}(G - w). \]
}
\begin{proof}
  Combining Proposition~\ref{S-eqs}, Theorem~\ref{S-swapc} and Corollary~\ref{S-bij-cor} gives the result.
\end{proof}
\vspace{\baselineskip}

%% file: r-case.tex
\section{The \texorpdfstring{$R$}{R}-case}\label{S:R-case}

Let $G$ be a connected 4-regular graph, and let $v$ and $w$ be two adjacent vertices of $G$ such that they do not share any neighbours.
Let $R = G - \{v,w\}$ be the graph obtained by removing vertices $v$ and $w$ from $G$, the last grey blob in Figure~\ref{fig:TSR}, with the neighbours of $v$ and $w$ labelled $\{a,b,c,d,e,f\}$ as in the figure.
Here $w$ has neighbours $\{a,b,c,v\}$, and $v$ has neighbours $\{d,e,f,w\}$.

The difficulty of the $R$-case stems from the fact that $v$ and $w$ no longer share any neighbours, and thus there is no distinguishing feature of $R$ that we can readily exploit.
However, the lack of specialness of the vertices $\{a,b,c,d,e,f\}$ is what lends itself to the symmetric nature of the $R$-case.
To prove the conjecture for $p=2$, we needed to use this symmetry to our advantage.

As we no longer have any 2-valent vertex to swap around, we look towards the more general control-vertex argument.
Notice in the swapping around a control vertex argument for the $S$-case there was a symmetric flavour in how the new edge partitions that we ended up with could have been from any $\mathcal{S}$ set of a similar form to the original, but we didn't need to know exactly which one.
Using a simplified version of the involution in the $S$-case, we use this symmetry to partially deal with the $R$-case.
For the rest of the sets, a more complicated control vertex argument is necessary, this time involving multiple control vertices.

\bpoint{Set-up}

As with the previous cases, we start with defining the sets of edge partitions that we are counting.

\tpointn{Definition} (Definition 3.1 of~\cite{specialc2})\label{R-defs}
\statement{
  Suppose $P$ is a bipartition of $\{a,b,c,d,e,f\}$. \\
  Let $\mathcal{R}_P$ be the set of bipartitions $(\psi, \phi)$ of the edges of $R$ such that $\psi$ is a spanning tree and $\phi$ is a spanning 2-forest compatible with $P$.
  Let $r_P = \abs{\mathcal{R}_P}$.
}

This time because of the symmetry of the neighbours of $v$ and $w$, and since we are counting modulo $2$, instead of picking specific vertex bipartitions for the spanning 2-forests we can add them all together.

\tpointn{Proposition} (Proposition 3.2, 3.3 of~\cite{specialc2})\label{R-eqs}
\statement[eq]{
  When $v$ and $w$ have no common neighbours
  \[ c_2^{(2)}(G - v) = r_{\{a\},\{b,c\}} + r_{\{b\},\{a,c\}} + r_{\{c\},\{a,b\}} \mod 2, \]
  \[ c_2^{(2)}(G - w) = r_{\{d\},\{e,f\}} + r_{\{e\},\{d,f\}} + r_{\{f\},\{d,e\}} \mod 2, \]
  and thus we have
  \begin{align*}\refstepcounter{equation}\tag{\theequation}\label{eq:R-counts}
    c_2^{(2)}(G - v) - c_2^{(2)}(G - w)
      =\ &{r}_{\{a\},\{b,c,d,e,f\}} + {r}_{\{b,c\},\{a,d,e,f\}} \\
      &+ {r}_{\{b\},\{a,c,d,e,f\}} + {r}_{\{a,c\},\{b,d,e,f\}} \\
      &+ {r}_{\{c\},\{a,b,d,e,f\}} + {r}_{\{a,b\},\{c,d,e,f\}} \\
      &+ {r}_{\{d\},\{a,b,c,e,f\}} + {r}_{\{e,f\},\{a,b,c,d\}} \\
      &+ {r}_{\{e\},\{a,b,c,d,f\}} + {r}_{\{d,f\},\{a,b,c,e\}} \\
      &+ {r}_{\{f\},\{a,b,c,d,e\}} + {r}_{\{d,e\},\{a,b,c,f\}} \mod 2. \\
  \end{align*}
}

\bpoint{Swapping around a control vertex}

Immediately from Equation~\eqref{eq:R-counts} we see that there are 6 terms corresponding to sets of the form $\mathcal{R}_{\{y,z\},\{x,*,*,*\}}$ where $\{x,y,z\}$ is either $\{a,b,c\}$ or $\{d,e,f\}$.
The thing to note is that these vertex bipartitions for the 2-forests are exactly in the form where we can apply Lemma~\ref{control-vertex}, this time using the other specification for the special vertex $x$.
The proof of this theorem uses a simpler version of the control vertex argument.

\tpointn{Theorem} (Lemma 4.4 of~\cite{specialc2})\label{R-control-bij}
\statement[eq]{
  There is a fixed-point free involution on
  \begin{align*}
    \mathcal{R}_{\{a,b\},\{c,d,e,f\}} &\cup \mathcal{R}_{\{a,c\},\{b,d,e,f\}} \cup \mathcal{R}_{\{b,c\},\{b,d,e,f\}} \\
    &\cup \mathcal{R}_{\{d,e\},\{a,b,c,f\}} \cup \mathcal{R}_{\{d,f\},\{a,b,c,e\}} \cup \mathcal{R}_{\{e,f\},\{a,b,c,d\}},
  \end{align*}
  and thus we have
  \begin{align*}
    {r}_{\{a,b\},\{c,d,e,f\}} &+ {r}_{\{a,c\},\{b,d,e,f\}} + {r}_{\{b,c\},\{b,d,e,f\}} \\
    &+ {r}_{\{d,e\},\{a,b,c,f\}} + {r}_{\{d,f\},\{a,b,c,e\}} + {r}_{\{e,f\},\{a,b,c,d\}} \equiv 0 \mod 2.
  \end{align*}
}

\bpoint{Two control vertices}\label{SS:R-new}

For the remaining 6 sets in Equation~\eqref{eq:R-counts}, we are no longer in a situation where we can apply Lemma~\ref{control-vertex} since these sets are of the form $\mathcal{R}_{\{x\},\{*,*,*,*,*\}}$ for $x \in \{a,b,c,d,e,f\}$.
What we need is a new notion of "control vertex", this time for a tree with five marked vertices instead of four.
However, as we no longer have any distinguishable vertices in this 5-part, we could not appropriately define one unique vertex as the control vertex.
Rather two control vertices are necessary to get an involution on the correct sets.

\begin{figure}[t]
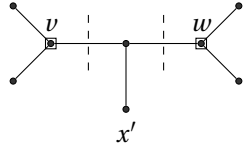

  \centering
  \ShapeFiveCut
  \caption[The two control vertices for Lemma~\ref{R-control}.]{The two control vertices $v$ and $w$ for Lemma~\ref{R-control}. The unlabelled leaves and $x'$ are all vertices in $p$. The dashed lines indicate the edges we will be picking for one of the involutions in the $R$-case. Some of the leaf edges are possibly contracted, with the exception that the edges of two leaves adjacent to the same vertex cannot simultaneously both be contracted.}
  \label{fig:two-control}
\end{figure}

\tpointn{Lemma}\label{R-control}
\statement{
  Suppose we have a bipartition of $\{a,b,c,d,e,f\}$ of the form $\{x\}, \{*,*,*, *, * \}$.
  Let $(\psi,\phi)$ be an edge partition in $\mathcal{R}_{\{x\},\{*,*,*,*,*\}}$ where $\psi$ is a spanning tree and $\phi$ is a spanning 2-forest compatible with $\{x\},\{*,*,*,*,*\}$.
  Let $t$ be the tree corresponding to part $p = \{*,*, *, *, *\}$ in the 2-forest $\phi$. \\\\
  Then, there is a unique pair of non-adjacent vertices $v$ and $w$ called the \textbf{control vertices}, such that removing $v$ and $w$ from $t$ partitions $p$ into all singletons.
  Figure~\ref{fig:two-control} depicts the two control vertices.\\
  Furthermore, in $t$ each control vertex will either be
  \begin{itemize}
    \item 2-valent and in $\{a,b,c,d,e,f\} \setminus \{x\}$, or
    \item 3-valent and not a vertex in $\{a,b,c,d,e,f\}$.
  \end{itemize}
  In all cases, the control vertices are leaves in the spanning tree $\psi$.
}
\begin{proof}
  In a similar fashion to the proof of Lemma~\ref{control-vertex}, we start by looking at the tree $t$ and the subtree created by taking the union of all the paths in $t$ between any two vertices in $p$.
  Since $p$ has five vertices this subtree will be of the form
  \[ \ShapeFive \]

  where the leaves are vertices in $p$, the edges are paths in $t$, and some of the paths may be contracted.

  Since $\psi$ is a spanning tree, the vertices in $\phi$, and thus $t$, are at most 3-valent, and the vertices in $\{a,b,c,d,e,f\}$, and thus $p$, are at most 2-valent.
  In particular, we have that the above is the only possible form for the subtree, where some of the paths may be contracted.
  Note trivially, the paths between two leaves adjacent to the same non-leaf vertex cannot be contracted simultaneously.
  Now, the valency restrictions also mean that in the subtree form above, the paths from $v$ and $w$ to the middle non-leaf vertex cannot be contracted.
  Thus $v \neq w$, and furthermore, they are not adjacent in $\phi$ i.e. there must be at least two edges in the path from $v$ to $w$.
  By the valency of the vertices in $\phi$, $v$ and $w$ can only be 2-valent if they are in the set $\{a,b,c,d,e,f\}$ or otherwise 3-valent, making them both leaves in the spanning tree $\psi$.
  As $R$ has more than two vertices, $v$ and $w$ are also not adjacent in $\psi$, and thus in $R$ the control vertices are non-adjacent.

  Finally, removing $v$ and $w$ from $t$ clearly partitions $p$ into all singletons and these are the unique pair of vertices which do so.
  The number of components of $t - \{v,w\}$ will depend on whether $v$ or $w$ are in the set $\{a,b,c,d,e,f\}$.
\end{proof}

Now that we've established the notion of control vertices, we can use a swapping argument once more for one final fixed-point free involution to cover the remaining sets in Equation~\eqref{eq:R-counts}.
As there are two vertices to swap around, like in Theorem~\ref{S-bijection}, we need a more complicated involution than previously.

\tpointn{Theorem}\label{R-bijection}
\statement[eq]{
  There is a fixed-point free involution on
  \begin{align*}
    \mathcal{R}_{\{a\},\{b,c,d,e,f\}} &\cup \mathcal{R}_{\{b\},\{a,c,d,e,f\}} \cup \mathcal{R}_{\{c\},\{a,b,d,e,f\}} \\
    &\cup \mathcal{R}_{\{d\},\{a,b,c,e,f\}} \cup \mathcal{R}_{\{e\},\{a,b,c,d,f\}} \cup \mathcal{R}_{\{f\},\{a,b,c,d,e\}},
  \end{align*}
  and thus we have
  \begin{align*}
    {r}_{\{a\},\{b,c,d,e,f\}} &+ {r}_{\{b\},\{a,c,d,e,f\}} + {r}_{\{c\},\{a,b,d,e,f\}} \\
    &+ {r}_{\{d\},\{a,b,c,e,f\}} + {r}_{\{e\},\{a,b,c,d,f\}} + {r}_{\{f\},\{a,b,c,d,e\}} \equiv 0 \mod 2.
  \end{align*}
}
\begin{proof}
  By construction the edge partitions in the union of sets of the form $\mathcal{R}_{\{x\},\{*,*,*,*,*\}}$, where $x \in \{a,b,c,d,e,f\}$, and the $*$'s are the rest of the vertices $\{a,b,c,d,e,f\} \setminus \{x\}$, satisfy the conditions of Lemma~\ref{R-control}.
  Let $(\psi, \phi)$ be an edge partition in any set of the form $\mathcal{R}_{\{x\},\{*,*,*,*,*\}}$, and let $v,w$ be its control vertices.
  Let $t$ be the tree in $\phi$ corresponding to part $\{*,*,*,*,*\}$ which contains $v$ and $w$.
  Since the control vertices are leaves in $\psi$, let $\eta_v$ and $\eta_w$ be the edges incident to $v$ and $w$ in $\psi$, respectively, and let $n_v$ and $n_w$ be their neighbours.

  Consider the following involution giving a new edge partition $(\psi',\phi')$:
  \begin{enumerate}[label=(\arabic*)]
    \item\label{R-step-1} \emph{Swapping either stays in} -- If $n_v \in t$ or $n_w \in t$:
      \begin{enumerate}[label=(\roman*), ref=(1\roman*)]
        \item\label{R-step-vin} \emph{$v$ stays in} -- If $n_v \in t$, so control vertex $v$ and its neighbour in $\psi$ are in the same tree of $\phi$, let $\eta_1$ be the edge incident to $v$ in $t$ in the path to $n_v$.\\
          Then, swap the edges $\eta_v$ and $\eta_1$ between $\psi$ and $\phi$.
        \item\label{R-step-win} \emph{$w$ stays in} -- If $n_w \in t$, so control vertex $w$ and its neighbour in $\psi$ are in the same tree of $\phi$, let $\eta_2$ be the edge incident to $w$ in $t$ in the path to $n_w$.\\
          Then, swap the edges $\eta_w$ and $\eta_2$ between $\psi$ and $\phi$.
      \end{enumerate}
      Here we are using non-exclusive or, so if $n_v, n_w \in t$, then swap around both $v$ and $w$ as above.
    \item\label{R-step-2} \emph{Swapping both goes out} -- Otherwise we have $n_v, n_w \not\in t$, where the control vertices and their neighbours in $\psi$ are in different trees of $\phi$.
      Let $\eta_1$ be the edge incident to $v$ in $t$ in the path to $w$, and let $\eta_2$ be the edge incident to $w$ in $t$ in the path to $v$; in Figure~\ref{fig:two-control} these are the edges indicated by the dashed lines. \\
      Then, swap edges $\eta_v$ and $\eta_w$ with $\eta_1$ and $\eta_2$ between $\psi$ and $\phi$.
  \end{enumerate}

  Firstly, as we are always changing which edges are incident to $v$ and/or $w$, $(\psi', \phi')$ must be different from $(\psi,\phi)$.
  Thus this transformation is fixed-point free.

  For~\ref{R-step-1}, both cases are symmetric so without loss of generality assume $n_v \in t$.
  Then, removing edge $\eta_1$ from $\phi$ breaks $t$ into two components and because of how $\eta_1$ was chosen, $n_v$ is now in a different component than $v$.
  As $v$ and its neighbour in $\psi$, $n_v$, are in the same tree of $\phi$, adding edge $\eta_v$ reconnects the two components from $t$.
  By the choice of $\eta_1$, we cannot create any cycles by adding $\eta_v$, and thus $\phi'$ is a spanning 2-forest which is compatible with the same vertex bipartition as $\phi$.
  Since $v$ is a leaf in $\psi$, removing edge $\eta_v$ and adding edge $\eta_1$ to create $\psi'$ maintains the spanning tree structure.
  Thus $(\psi', \phi')$ is a valid edge partition in the same $\mathcal{R}$ set as $(\psi, \phi)$.
  We also have if the control vertices stay the same, $(\psi', \phi')$ will still have the property that $v$ and its neighbour in $\psi'$ are in the same tree of $\phi'$, and thus applying the transformation again returns to $(\psi, \phi)$.

  When both $n_v, n_w \in t$, since $v\neq w$ and they are not adjacent, the four edges in the swap are all distinct; $\eta_1 \neq \eta_2$ and $\eta_v \neq \eta_w$.
  Thus we can simultaneously swap edges $\eta_v$ and $\eta_w$ with edges $\eta_1$ and $\eta_2$, respectively.
  Note that we need to swap around both control vertices as there is no way to distinguish $v$ from $w$, and thus no way to ensure we get an involution if we arbitrarily picked one to swap edges around.
  Applying the above swapping argument twice, we have that $(\psi', \phi')$ is a valid edge partition in the same $\mathcal{R}$ set as before and applying the transformation again returns to $(\psi, \phi)$.

  For~\ref{R-step-2} when $n_v, n_w \not\in t$, once again we know that the four edges in the swap are all distinct by the same reasoning as when $n_v, n_w \in t$.
  Then, removing edges $\eta_1$ and $\eta_2$ breaks $t$ into three subtrees (thus breaking $\phi$ into four components); one with $v$, one with $w$, and one with the vertex $x' \in \{a,b,c,d,e,f\} \setminus \{x\}$ as in Figure~\ref{fig:two-control}.
  Now since $v$ and $n_v$ are in different trees of $\phi$, adding $\eta_v$ connects the subtree from $t - \{\eta_1, \eta_2\}$ with $v$ to the tree in $\phi$ corresponding to part $\{x\}$.
  Similarily, as $w$ and $n_w$ are in different trees of $\phi$, adding $\eta_w$ connects the subtree with $w$ also to the tree with $x$ in $\phi$.

  As the last subtree, which contains $x'$, never gets connected to any other component, this transformation creates a new spanning 2-forest $\phi'$ which is compatible with the vertex bipartition $\{x'\}, \{x, *,*,*,*\}$.
  Since $v,w$ are leaves in $\psi$, removing edges $\eta_v$ and $\eta_w$ from $\psi$ disconnects the control vertices and adding edges $\eta_1$ and $\eta_2$ reconnects them, creating $\psi'$ which is once again a spanning tree.
  Thus $(\psi', \phi')$ is a valid edge partition in $\mathcal{R}_{\{x'\},\{x,*,*,*,*\}}$, which is in the correct form.
  Once again, if the control vertices stay the same, then their neighbours in $\psi$ will be in a different tree of $\phi'$ landing us back in ~\ref{R-step-2}.
  Thus applying the transformation again brings us back to $(\psi, \phi)$.

  What's left of the proof is to make sure $v$ and $w$ are still the control vertices in $(\psi', \phi')$.
  Notice we can formulate the control vertices as the two vertices in $t$ such that for a specific $x'$ in part $\{*,*,*,*,*\}$, each vertex is the last common vertex on the paths from $x'$ to the exactly two of the other vertices in the part.
  The vertex $x'$ is as depicted in Figure~\ref{fig:two-control}.
  Under this formulation, we immediately see that in~\ref{R-step-2} as $\phi'$ is created by connecting subtrees with $v$ and $w$ to the tree corresponding to part $\{x\}$ in $\phi$, the vertex $x$ now acts as the $x'$ for $\phi'$.
  Thus $v$ and $w$ remain the control vertices for $(\psi',\phi')$.

  For~\ref{R-step-1}, depending on where $n_v$ or $n_w$ were in the tree $t$, the vertex acting as $x'$ in $\phi'$ may no longer be the same $x'$ from $\phi$.
  However, one can check that in all possible situations, the property for the control vertices is still satisfied by $v$ and $w$, and thus $v$ and $w$ are still the control vertices.

  Finally putting everything together, we indeed get a fixed-point free involution on the union of sets of the form $\mathcal{R}_{\{x\},\{*,*,*,*,*\}}$,
  \begin{align*}
    \mathcal{R}_{\{a\},\{b,c,d,e,f\}} &\cup \mathcal{R}_{\{b\},\{a,c,d,e,f\}} \cup \mathcal{R}_{\{c\},\{a,b,d,e,f\}} \\
    &\cup \mathcal{R}_{\{d\},\{a,b,c,e,f\}} \cup \mathcal{R}_{\{e\},\{a,b,c,d,f\}} \cup \mathcal{R}_{\{f\},\{a,b,c,d,e\}},
  \end{align*}
  and therefore, the size of this must be even
  \begin{align*}
    {r}_{\{a\},\{b,c,d,e,f\}} &+ {r}_{\{b\},\{a,c,d,e,f\}} + {r}_{\{c\},\{a,b,d,e,f\}} \\
    &+ {r}_{\{d\},\{a,b,c,e,f\}} + {r}_{\{e\},\{a,b,c,d,f\}} + {r}_{\{f\},\{a,b,c,d,e\}} \equiv 0 \mod 2.
  \end{align*}
\end{proof}

Notice that we can think of this involution involving two control vertices as a generalized version of the one used in Theorem~\ref{S-bijection}.
For that involution, the 2-valent vertex $c$ acts as a second "control vertex".

\bpoint{Completing the \texorpdfstring{$R$}{R}-case}

While there were no special vertices in $\{a,b,c,d,e,f\}$ that we could readily use to swap between specific vertex bipartitions, we were able to take advantage of the symmetric nature of the $R$-case and use control vertices to find involutions where we only cared about the form of the vertex bipartition we may be swapping our edge partitions to.
 With the two swapping arguments, one involving one control vertex and the other involving two control vertices, we have everything we need to prove the $c_2$ completion conjecture for $p=2$ in the $R$-case.

\tpointn{Theorem}\label{R-case}
\statement[eq]{
  Let $G$ be a connected 4-regular graph.
  Let $v$ and $w$ be adjacent vertices of $G$ such that they do not share any common neighbours.
  Then,
  \[ c_2^{(2)}(G - v) = c_2^{(2)}(G - w). \]
}
\begin{proof}
  Together Proposition~\ref{R-eqs}, Theorem~\ref{R-control-bij} and Theorem~\ref{R-bijection} gives the result.
\end{proof}
\vspace{\baselineskip}

%% file: completing.tex
\section{Completing the \texorpdfstring{$p=2$}{p=2} case}\label{S:completing-conjecture}

To prove the $c_2$ completion conjecture for $p=2$, in all three cases the underlying idea remained the same: finding a vertex, or vertices, and picking particular edges incident to them to swap between the two parts of an edge partition, where this action gives a fixed-point free involution on specific sets of edge partitions.

In the $T$-case, we found that we could use the two 2-valent vertices to swap edges around (see ~\cite{Hmmath} specificially for the $p=2$ case, or see Section~\ref{S:higherp} for the general $p$ argument).
In the $S$-case, we reused the 2-valent vertex swapping argument and generalized it to give the notion of swapping around a particular vertex called the control vertex.
We found that for a second involution, we needed a two-phase swapping process where we used both the control vertex and a 2-valent vertex.
Finally in the $R$-case, we reused the control vertex argument and found that we needed to extend it to include two control vertices for a second involution, this extension itself being a generalized version of the second involution in the $S$-case.

We note here that the two new bijections, giving Theorems~\ref{S-bijection} and~\ref{R-bijection} in the $S$ and $R$ cases, respectively, were implemented in \texttt{Sage}~\cite{sagemath} to verify our results on some small graphs.

Assembling everything together, we finally complete the $c_2$ completion conjecture for $p=2$!

\tpointn{Theorem}\label{p2-completed}
\statement[eq]{
  Let $G$ be a connected 4-regular graph, and let $v$ and $w$ be vertices of $G$.
  Then,
  \[ c_2^{(2)}(G - v) = c_2^{(2)}(G - w). \]
}
\begin{proof}
  As in the discussion in Section~\ref{S:intro}, since $G$ is connected there is path between any two vertices of $G$, and thus it suffices to prove the conjecture when $v$ and $w$ are adjacent vertices.
  For any two non-adjacent vertices, we could then just follow any path between them, getting a chain of equivalences to obtain the required result.

  Then, as $G$ is a 4-regular graph, the four possibilities for the neighbours of $v$ and $w$ are:
  \begin{itemize}
    \item $v$ and $w$ share all neighbours, or
    \item $v$ and $w$ have exactly two common neighbours -- $T$-case, or
    \item $v$ and $w$ have only one common neighbour -- $S$-case, or
    \item $v$ and $w$ do not share any neighbours -- $R$-case.
  \end{itemize}

  When $v$ and $w$ share all neighbours, notice that $G - v$ and $G - w$ are isomorphic.
  Thus, trivially, their $c_2$'s are equivalent.
  When $v$ and $w$ share some neighbours or none at all, Corollary~\ref{T-case}, Theorem ~\ref{S-case}, and Theorem~\ref{R-case} give the result in the $T$, $S$, and $R$ cases, respectively.
  This completes the proof.
\end{proof}
\vspace{\baselineskip}

%% file: higherp.tex
\section{Higher \texorpdfstring{$p$}{p} in the \texorpdfstring{$T$}{T}-case}\label{S:higherp}

The hope in finding these counting methods is, of course, to be able to generalize the results in the $p=2$ case to higher values of $p$, which would then prove the $c_2$ completion conjecture in full.
Optimistically, we do know how to completely finish the $T$-case for all $p$ and this same argument gives some partial results in the $S$-case.
While the (full) $S$ and $R$ cases are not as simple, this extension gives some hope for generalizing the arguments needed in those cases for higher values of $p$.
Note that the double triangle argument, see Section \ref{S:T-case}, for the $T$-case does still apply, but since it only reduces the question to a smaller graph, and we do not have the full result on all graphs, a direct proof for the $T$-case is stronger.
We present the argument here.

\bpoint{Set-up}

Recall that in the $T$-case, $T = G - \{v,w\}$ where adjacent vertices $v$ and $w$ share two common neighbours labelled $b$ and $c$ (see Figure~\ref{fig:TSR}).
Thus in $T$, vertices $b$ and $c$ are now 2-valent.

From Equation~\eqref{eq:c2-higher} and the discussion thereafter, for $p > 2$ we are now counting the number of ways to partition $p-1$ copies of each edge of $G - \{v,w\}$ into $2p-2$ parts, where the first $p-1$ parts are spanning trees and the last $p-1$ parts are spanning 2-forests compatible with some partition of the neighbours of the vertex we reduced by.
The difficulty of the higher $p$ case comes from the fact that, unlike the $p=2$ case where we had graph-complement pairs $(\psi, \phi)$ where $\phi = \overline{\psi}$, now we can have spanning 2-forests or spanning trees that did not previously appear in the $p=2$ case.

\tpointn{Definition} \label{T-higher-defs}
\statement{
  Let $p$ be any prime.
  Suppose $P_1,\ldots, P_{p-1}$ are bipartitions of $\{a,b,c,d\}$. \\
  Let $T^{(p-1)}$ denote the graph with $(p-1)$-copies of each edge in $T$.
  Let $\mathcal{T}_{P_1,\ldots,P_{p-1}}$ be the set of (ordered) partitions $(\psi_1,\ldots,\psi_{p-1}, \phi_1,\ldots, \phi_{p-1})$ of the edges of $T^{(p-1)}$ such that for each $i$, $\psi_i$ is a spanning tree in $T$ and $\phi_i$ is a spanning 2-forest in $T$ compatible with $P_i$.
  Let $t_{P_1,\ldots,P_{p-1}} = \abs{\mathcal{T}_{P_1,\ldots,P_{p-1}}}$.
}

Like in the $p=2$ case, we can choose the order of the edges incident to the 3-valent vertex we are reducing by in Equation~\eqref{eq:c2-higher} to exploit that $v$ and $w$ share two common neighbours.

\tpointn{Proposition}\label{T-higher-eqs}
\statement[eq]{
  Let $p$ be any prime.
  Let $P = \{a\},\{b,c,d\}$, $P' = \{d\},\{a,b,c\}$ and $Q = \{a,d\},\{b,c\}$.\\
  When $v$ and $w$ have common neighbours $b$ and $c$,
  \[ c_2^{(p)}(G-v)
      = -\sum_{P_i = P \text{ or } Q} t_{P_1,\ldots,P_{p-1}}
      = -\sum_{\ell=0}^{p-1} \binom{p-1}{\ell} t_{\scriptsize\underbrace{P,\ldots,P}_{\ell},\underbrace{Q,\ldots,Q}_{p-1-\ell}} \mod p, \]
  \[ c_2^{(p)}(G-w)
      = -\sum_{P_i = P' \text{ or } Q} t_{P_1,\ldots,P_{p-1}}
      = -\sum_{\ell=0}^{p-1} \binom{p-1}{\ell} t_{\scriptsize\underbrace{P',\ldots,P'}_{\ell},\underbrace{Q,\ldots,Q}_{p-1-\ell}} \mod p, \]
  and thus we have
  \[ c_2^{(p)}(G - w) - c_2^{(p)}(G - v) = \sum_{\ell=1}^{p-1} \binom{p-1}{\ell} \left( t_{\scriptsize\underbrace{P,\ldots,P}_{\ell},\underbrace{Q,\ldots,Q}_{p-1-\ell}} - t_{\scriptsize\underbrace{P',\ldots,P'}_{\ell},\underbrace{Q,\ldots,Q}_{p-1-\ell}} \right) \mod p. \]
}
\begin{proof}
Reducing with respect to $w$ for $G - v$, Equation~\eqref{eq:c2-higher} tells us
\begin{align*}
  c_2^{(p)}(G-v)
    &\equiv -[\alpha_4^{p-1}\cdots\alpha_{|E|}^{p-1}]\; \left(\Phi_{T}^{\{a\},\{b,c\}} \Psi_{T} \right)^{p-1} \mod p \\
    &\equiv -[\alpha_4^{p-1}\cdots\alpha_{|E|}^{p-1}]\; \left(\Phi_{T}^{\{a\},\{b,c,d\}} + \Phi_T^{\{a,d\},\{b,c\}}\right)^{p-1}\Psi_{T}^{p-1} \mod p \\
    &\equiv -\sum_{P_i = P \text{ or } Q} t_{P_1,\ldots,P_{p-1}} \mod p
\end{align*}
where the sum in the last equation runs over all $2^{p-1}$ possible tuples $(P_1,\ldots,P_n)$ where $P_i$ is one of the two partitions $P$ or $Q$.

Equivalently, we can write
\begin{align*}
  c_2^{(p)}(G-v)
    &\equiv -[\alpha_4^{p-1}\cdots\alpha_{|E|}^{p-1}]\; \sum_{\ell=0}^{p-1} \binom{p-1}{\ell} \left(\Phi_{T}^{\{a\},\{b,c,d\}}\right)^{\ell}\left(\Phi_T^{\{a,d\},\{b,c\}}\right)^{p-1-\ell} \Psi_{T}^{p-1} \mod p \\
    &\equiv -\sum_{\ell=0}^{p-1} \binom{p-1}{\ell} t_{\scriptsize\underbrace{P,\ldots,P}_{\ell},\underbrace{Q,\ldots,Q}_{p-1-\ell}} \mod p
\end{align*}
where now we are fixing the tuples $(P_1,\ldots,P_{p-1})$ to those where the first $\ell$ bipartitions are $P$ and the last $p-1-\ell$ bipartitions are $Q$.

Similarily, reducing with respect to $v$ for $G - w$ gives
\begin{align*}
  c_2^{(p)}(G-w)
    &\equiv -[\alpha_4^{p-1}\cdots\alpha_{|E|}^{p-1}]\; \left(\Phi_{T}^{\{d\},\{b,c\}} \Psi_{T} \right)^{p-1} \mod p \\
    &\equiv -[\alpha_4^{p-1}\cdots\alpha_{|E|}^{p-1}]\; \left(\Phi_{T}^{\{d\},\{a,b,c\}} + \Phi_T^{\{a,d\},\{b,c\}}\right)^{p-1}\Psi_{T}^{p-1} \mod p \\
    &\equiv -\sum_{P_i = P' \text{ or } Q} t_{P_1,\ldots,P_{p-1}} \mod p
\end{align*}
where the sum in the last equation runs over all $2^{p-1}$ possible tuples $(P_1,\ldots,P_n)$ where $P_i$ is now one of the two partitions $P'$ or $Q$.
Equivalently, we can write
\begin{align*}
  c_2^{(p)}(G-w)
    &\equiv -\sum_{\ell=0}^{p-1} \binom{p-1}{\ell} t_{\scriptsize\underbrace{P',\ldots,P'}_{\ell},\underbrace{Q,\ldots,Q}_{p-1-\ell}} \mod p.
\end{align*}

Then subtracting the two gives
\begin{align*}
  c_2^{(p)}(G - w) - c_2^{(p)}(G - v)
    &\equiv \sum_{\ell=0}^{p-1} \binom{p-1}{\ell} t_{P,\ldots,P,Q,\ldots,Q} - \sum_{\ell=0}^{p-1} \binom{p-1}{\ell} t_{P',\ldots,P',Q,\ldots,Q} \mod p \\
    &\equiv \sum_{\ell=1}^{p-1} \binom{p-1}{\ell} \left( t_{\scriptsize\underbrace{P,\ldots,P}_{\ell},\underbrace{Q,\ldots,Q}_{p-1-\ell}} - t_{\scriptsize\underbrace{P',\ldots,P'}_{\ell},\underbrace{Q,\ldots,Q}_{p-1-\ell}} \right) \mod p
\end{align*}
where we noticed that when $\ell = 0$, we have $t_{Q,\ldots,Q}$ in both sums.
\end{proof}

\bpoint{Swapping around two-valent vertices}

The argument proceeds by swapping edges around the two 2-valent vertices, analogous to swapping around $c$ in the $S$-case.

\tpointn{Lemma}\label{higher-swap-two}
\statement{
  Let $V$ be a set of marked vertices of a connected graph $G$, and suppose $c \in V$ is 2-valent.
  Let $G^{(p-1)}$ be the graph with $(p-1)$-copies of each edge in $G$.
  Suppose $\pi = (\psi_1,\ldots, \psi_{p-1},\phi_1,\ldots,\phi_{p-1})$ is an edge partition of $G^{(p-1)}$ where each $\psi_i$ is a spanning tree in $G$ and each $\phi_i$ is a spanning 2-forest in $G$ compatible with a bipartition of $V$ where every part contains at least one of the vertices in $V \setminus \{c\}$. \\\\
  Then, of the the two edges incident to $c$ in $G$, exactly one copy of the two edges is in each part of $\pi$.
  Furthermore, if we pick a spanning 2-forest $\phi_i$ and a spanning tree $\psi_j$ such that the two have different edges incident to $c$, then swapping these two edges between $\phi_i$ and $\psi_j$ yields a new edge partition $\pi'$ with the same properties, except possibly the bipartition of $V$ that $\phi'_i$ is compatible with.
  Note that such a spanning tree $\psi_j$ always exists.
}
\begin{proof}
  Looking at the partition $\pi$, there must be at least one edge incident to $c$ in each $\psi_i$ as they are all spanning trees in $G$.
  Now as each $\phi_i$ is a spanning 2-forest in $G$ compatible with a bipartition of $V$ where $c$ is in a part with at least one other vertex in $V$, $c$ cannot be an isolated vertex in any $\phi_i$.
  Thus each $\phi_i$ must have at least one edge incident to $c$.
  As $c$ is 2-valent in $G$, there are only $2(p-1)$ edges incident to $c$ in $G^{(p-1)}$ to be distributed among the $2(p-1)$ parts of $\pi$.
  Since each of these parts must have at least one edge incident to $c$, we have that exactly one edge is incident to $c$ in each $\psi_i$ and each $\phi_i$.
  That is, $c$ is a leaf in every part of $\pi$.
  Furthermore, as there are only $p-1$ copies of each edge in $G^{(p-1)}$, there is at least one spanning tree $\psi_j$ such that $\psi_j$ has a different edge incident to $c$ as $\phi_i$.

  For the last part, we pick a $\psi_j$ and a $\phi_i$ such that they have different edges incident to $c$ and fix the rest of the parts.  Note that such a $\psi_j$ must always exist because there are only $p-1$ copies of each edge to distribute.
  Because $c$ is a leaf in both $\psi_i$ and $\phi_j$, removing the edges incident to $c$ disconnects and furthermore, isolates $c$ from the rest of the tree and the 2-forest, respectively.
  Since both $\psi_j$ and $\phi_i$ are spanning, they must contain both neighbours of $c$, and thus we can reconnect $c$ via the opposite edges incident to $c$ as in $\psi_j$ and $\phi_i$.
  This gives a new edge partition $\pi'$ which only differs from $\pi$ at the $j$-th spanning tree $\psi'_j$ and $i$-th spanning 2-forest $\phi'_i$, which are $\psi_j$ and $\phi_i$ respectively but with the edges incident to $c$ swapped.
  However, in $\phi'_i$, $c$ could be reconnected to a different tree than in $\phi_i$ thus changing the bipartition of $V$ that $\phi'_i$ is compatible with.
  Specifically, either the bipartition of $V$ remains the same as for $\phi_i$ or $c$ swaps between the parts of the bipartition.
\end{proof}

For a general connected graph $G$, we analogously define the $\mathcal{G}_{P_1,\ldots,P_{p-1}}$ and its size $g_{P_1,\ldots,P_{p-1}}$ as in Definition~\ref{T-higher-defs}, where here each $P_i$ is a bipartition of some set of marked vertices $V$ of $G$.
With Lemma~\ref{higher-swap-two}, we can extend the swapping around a 2-valent vertex argument for $p \geq 2$.
The main observation is that we can pick one spanning 2-forest and count all the ways to swap around a 2-valent vertex with the $p-1$ available spanning trees, while fixing the rest of the spanning forests.

\tpointn{Theorem}\label{higher-swap}
\statement[eq]{
  Let $V$ be a set of marked vertices of a connected graph $G$ and suppose $c \in V$ is 2-valent.
  Suppose $P_1,\ldots,P_{p-1}$ are bipartitions of $V$ such that every part contains at least one of the vertices in $V \setminus \{c\}$.
  Then for any index $i$ there is a partition of the tuples in
  \[ \mathcal{G}_{P_1,\ldots,P_i, \ldots, P_{p-1}} \cup \mathcal{G}_{P_1,\ldots,P_i^{c},\ldots,P_{p-1}} \]
  into orbits whose sizes are divisible by $p$, where $P_i^c$ is the bipartition of $V$ obtained from $P_i$ by swapping which part $c$ is in.
  Note that $P_1,\ldots,P_{i-1},P_{i+1},\ldots,P_{p-1}$ remain fixed.
  Thus we have that
  \[ g_{P_1,\ldots,P_i, \ldots, P_{p-1}} + g_{P_1,\ldots,P_i^{c},\ldots,P_{p-1}} \equiv 0 \mod p. \]
}
\begin{proof}
  First, as $c$ is a 2-valent vertex in $T$, every partition in $\mathcal{G}_{P_1,\ldots,P_i, \ldots, P_{p-1}} \cup \mathcal{G}_{P_1,\ldots,P_i^{c},\ldots,P_{p-1}}$ satisfies the conditions in Lemma~\ref{higher-swap-two}.
  Consider any $\pi = (\psi_1,\ldots,\psi_{p-1}, \phi_1,\ldots,\phi_{p-1}) \in \mathcal{G}_{P_1,\ldots,P_i, \ldots, P_{p-1}}$ and fix some index $i$, choosing a spanning 2-forest $\phi_i$ which is compatible with $P_i$.
  By Lemma~\ref{higher-swap-two}, we have that $c$ is a leaf in every part of $\pi$ and for any spanning tree $\psi_j$ that has a different edge incident to $c$ than $\phi_i$, we can swap the edges incident to $c$ between $\phi_i$ and $\psi_j$ to get a new edge partition with the same properties.
  If the two neighbours of $c$ were in the same tree in $\phi_i$, any swap gives a different edge partition in $\mathcal{G}_{P_1,\ldots,P_i, \ldots, P_{p-1}}$.
  Otherwise, if the two neighbours were in different trees in $\phi_i$, we get an edge partition now in $\mathcal{G}_{P_1,\ldots,P_i^c, \ldots, P_{p-1}}$.
  In either case, swapping around $c$ again returns us to the original partition in $\mathcal{G}_{P_1,\ldots,P_i, \ldots, P_{p-1}}$.

  What's left is to count the size of the orbits created from this swapping action between $\phi_i$ and the $p-1$ available spanning trees.
  Let $\eta_1$ and $\eta_2$ be the two edges incident to $c$ with $\eta_1 \in \phi_i$.
  Suppose there are $k$ spanning trees in $\psi_1,\ldots,\psi_{p-1}$ that contain $\eta_2$ and thus $p-1-k$ spanning trees that contain $\eta_1$.
  As there are only $p-1$ copies of each edge distributed amongst the parts of $\pi$ and $\phi_i$ contains one copy of $\eta_1$, we know that $1 \leq k \leq p-1$.
  Picking any such $\psi_j$ that contains $\eta_2$, swapping $\eta_1$ and $\eta_2$ between $\phi_i$ and $\psi_j$ gives a new edge partition $\pi'$ which is $\pi$ except the $j$-th spanning tree $\psi'_j$ is $\psi_j$ with the edge $\eta_1$ instead of $\eta_2$ and the $i$-th spanning 2-forest $\phi'_i$ is $\phi_i$ with the edge $\eta_2$ swapped for $\eta_1$.

  Looking at $\pi'$, now $\phi'_i$ contains $\eta_2$ and there are $p-k$ spanning trees with $\eta_1$ (the same $p-1-k$ spanning trees as $\pi$ and also $\psi'_j$).
  Thus there are $p-k$ ways to swap around $c$ to get back to the spanning 2-forest $\phi_i$ in the $i$-th position, one of which returns us exactly to the edge partition $\pi$.
  We can then apply the swapping argument again to the new edge partitions.

  In particular, every partition in the orbit containing $\pi$ will have as spanning 2-forests either $\{\phi_1,\ldots,\phi_{p-1}\}$ or $\{\phi_1,\ldots,\phi'_1,\ldots,\phi_{p-1}\}$.
  The spanning trees of the partitions will be $\{\psi_1,\ldots,\psi_{p-1}\}$ except with possibly different edges incident to $c$ in each.
  Thus, to count the size of the orbit containing $\pi$, what we are really counting is the number of ways to distribute the $p-k$ copies of $\eta_1$ and $k$ copies of $\eta_2$ between the $i$-th spanning 2-forest and the $p-1$ spanning trees.
  As we are distributing $k$ copies of $\eta_2$ into $p$ slots, which then automatically fixes the $p-k$ slots with $\eta_1$, there are $\binom{p}{k}$ ways to do this distribution.
  Therefore the size of the orbit is exactly $\binom{p}{k}$ and as $1 \leq k \leq p-1$, we have that $p$ divides $\binom{p}{k}$.

  To conclude, we have shown is that we can partition $\mathcal{G}_{P_1,\ldots,P_i, \ldots, P_{p-1}} \cup \mathcal{G}_{P_1,\ldots,P_i^{c},\ldots,P_{p-1}}$, into orbits of size divisible by $p$ and thus
  \[ g_{P_1,\ldots,P_i, \ldots, P_{p-1}} + g_{P_1,\ldots,P_i^{c},\ldots,P_{p-1}} \equiv 0 \mod p. \]
\end{proof}

\tpointn{Remark}\label{higher-swap-rem}
\statement{
  We actually know more information about each of the orbits in
  \[ \mathcal{G}_{P_1,\ldots,P_i, \ldots, P_{p-1}} \cup \mathcal{G}_{P_1,\ldots,P_i^{c},\ldots,P_{p-1}}. \]
  Let $\pi = (\psi_1,\ldots,\psi_{p-1}, \phi_1,\ldots,\phi_{p-1})$ be any edge partition in this set and let $k$ be the number of spanning trees in $\pi$ that have a different edge incident to $c$ than $\phi_i$.
  \begin{itemize}[itemsep=4pt]
    \item The orbit of $\pi$ is of size $\binom{p}{k}$.
    \item If the two neighbours of $c$ are in the same tree in $\phi_i$, the orbit of $\pi$ is fully contained in $\mathcal{G}_{P_1,\ldots,P_i, \ldots, P_{p-1}}$.
      In this case, we say that $\phi_i$ "swaps in" with respect to $c$, and the orbits of all such $\pi$ form a set whose size is divisible by $p$ in $\mathcal{G}_{P_1,\ldots,P_i, \ldots, P_{p-1}}$.
    \item If the two neighbours of $c$ are in different trees in $\phi_i$, the orbit of $\pi$ is split across the two sets with exactly $\binom{p-1}{p-1-k}$ partitions in $\mathcal{G}_{P_1,\ldots,P_i, \ldots, P_{p-1}}$ and $\binom{p-1}{p-k}$ in $\mathcal{G}_{P_1,\ldots,P_i^c, \ldots, P_{p-1}}$.
  \end{itemize}
  In fact, we can relax the condition on the bipartitions $P_1,\ldots,P_{p-1}$ if we are considering the orbits of specific edge partitions $\pi$.
  All we need is for $c$ to be a leaf in every part of $\pi$.
  Then, the swapping argument still holds and we obtain the same results on the orbit of $\pi$ as above.
}

When $p=2$, note that Theorem~\ref{higher-swap} implies Theorem~\ref{S-swapc} in the $S$ case.

\tpointn{Corollary}\label{T-higher-swap}
\statement{
  Let $p$ be any prime.
  For any $\ell = 1, \ldots, p-1$, we have
  \[ t_{\scriptsize\underbrace{P,\ldots,P}_{\ell},\underbrace{Q,\ldots,Q}_{p-1-\ell}} \equiv t_{\scriptsize\underbrace{P',\ldots,P'}_{\ell},\underbrace{Q,\ldots,Q}_{p-1-\ell}} \mod p \]
  for partitions $P = \{a\},\{b,c,d\}$, $P' = \{d\},\{a,b,c\}$ and $Q = \{a,d\},\{b,c\}$.
}
\begin{proof}
  Recall that in the graph $T$, both $b$ and $c$ are 2-valent vertices.
  The key observation is that $P$, $P'$ and $Q$ are all bipartitions of $V = \{a,b,c,d\}$ that satisfy the condition in Theorem~\ref{higher-swap} for $b$ and for $c$.
  Furthermore, swapping which part $b$ is in in $P$ and swapping which part $c$ is in in $P'$ yields the same bipartition of $V$, $P^b = \{a,b\},\{c,d\} = P'^c$, where we are using the superscript to indicate which element changes part.  The bipartition $\{a,b\},\{c,d\}$ also satisfies the condition that every part contains at least one of the vertices in $V \setminus \{b\}$ and similarly for $V \setminus \{c\}$.
  Thus we can repeatedly apply Theorem~\ref{higher-swap} with $G = T$ and $V = \{a,b,c,d\}$ to both sides of the equation.

  Starting with $\{P_1,\ldots,P_{p-1}\} = \{P,\ldots,P,Q,\ldots,Q\}$ where there are $\ell \geq 1$ $P$'s, we apply Theorem~\ref{higher-swap} with the $P$ at index $\ell$.
  Swapping around $b$ in $P_{\ell} = P$ then gives
  \[ t_{\scriptsize\underbrace{P,\ldots,P}_{\ell},\underbrace{Q,\ldots,Q}_{p-1-\ell}} + t_{\scriptsize\underbrace{P,\ldots,P}_{\ell-1},P^b,\underbrace{Q,\ldots,Q}_{p-1-\ell}} \equiv 0 \mod p \]
  If $\ell - 1 > 0$, we now apply the theorem with $\{P_1,\ldots,P_{p-1}\} = \{P,\ldots,P,P^b,Q,\ldots,Q\}$ and swapping around $b$ in $P_{\ell-1} = P$, to get
  \[ t_{\scriptsize\underbrace{P,\ldots,P}_{\ell-1},P^b,\underbrace{Q,\ldots,Q}_{p-1-\ell}} + t_{\scriptsize\underbrace{P,\ldots,P}_{\ell-2},P^b,P^b,\underbrace{Q,\ldots,Q}_{p-1-\ell}} \equiv 0 \mod p \]
  Continuing successively in this fashion and swapping around the rest of the $P$'s, we obtain
  \[ t_{\scriptsize\underbrace{P,\ldots,P}_{\ell},\underbrace{Q,\ldots,Q}_{p-1-\ell}} \equiv (-1)^{\ell} t_{\scriptsize\underbrace{P^b,\ldots,P^b}_{\ell},\underbrace{Q,\ldots,Q}_{p-1-\ell}} \mod p \]
  Similarily, starting with $\{P_1,\ldots,P_{p-1}\} = \{P',\ldots,P',Q,\ldots,Q\}$, we successively apply Theorem~\ref{higher-swap}, now swapping around $c$ in each $P'$, to get
  \[ t_{\scriptsize\underbrace{P',\ldots,P'}_{\ell},\underbrace{Q,\ldots,Q}_{p-1-\ell}} \equiv (-1)^{\ell} t_{\scriptsize\underbrace{P'^c,\ldots,P'^c}_{\ell},\underbrace{Q,\ldots,Q}_{p-1-\ell}} \mod p \]
  As $P^b = \{a,b\},\{c,d\} = P'^c$,
  \[ t_{\scriptsize\underbrace{P,\ldots,P}_{\ell},\underbrace{Q,\ldots,Q}_{p-1-\ell}} \equiv t_{\scriptsize\underbrace{P',\ldots,P'}_{\ell},\underbrace{Q,\ldots,Q}_{p-1-\ell}} \mod p \]
\end{proof}

\bpoint{Completing the \texorpdfstring{$T$}{T}-case}

By putting together Proposition~\ref{T-higher-eqs} and Corollary~\ref{T-higher-swap}, we immediately obtain the full completion conjecture in the $T$-case!

\tpointn{Theorem}\label{T-higher}
\statement[eq]{
  Let $G$ be a connected 4-regular graph and $p$ be any prime.
  Let $v$ and $w$ be adjacent vertices of $G$ such that they share two common neighbours.
  Then,
  \[ c_2^{(p)}(G - v) = c_2^{(p)}(G - w). \]
  and thus we have
  \[ c_2(G - v) = c_2(G - w). \]
} \\

Note that these arguments do give us a bit of leverage on the $S$-case as well.  Suppose we have $(\psi_1,\ldots,\psi_{p-1}, \phi_1,\ldots,\phi_{p-1})$ in the $S$-case.  Using Remark \ref{higher-swap-rem}, if for some $\phi_i$ swapping around $c$ stays in the same vertex bipartition, then we can assign the edges around $c$ from $\phi_i$ and from all the spanning trees to obtain $\binom{p}{k}$ edge partitions all associated to the same vertex bipartitions, where $1\leq k \leq p-1$. Hence modulo $p$ these edge partitions do not contribute.  On the other hand, if none of the $\phi_i$ swap in, and if the vertex bipartitions are such that swapping out gives a bipartition that occurs from decompleting the other way, then by swapping all of the trees and forests around $c$ we find we can match together these two terms, one from each side of the desired completion equation.  This works for those vertex bipartitions which only required swapping around $c$ in the $p=2$ proof.  

While how to extend this to the rest of the $S$-case and to the $R$-case remains unclear, these results are suggestive that the edge-partitioning approach to completion is viable for $p>2$.  Though the arguments will be more intricate, there are many opportunities for configurations to gather in sets of size divisible by $p$ on account of appropriate binomial coefficients.  The main difficulty seems to be that while swapping around a 2-valent vertex remains similar to the $p=2$ case, we can no longer necessarily find appropriate control vertices within the graph.\\